\documentclass[a4paper,onecolumn, superscriptaddress,10pt,accepted=2026-03-03,issue=3, volume=8, shorttitle=papers]{compositionalityarticle}
 \pdfoutput=1
\usepackage[utf8]{inputenc}
\usepackage[english]{babel}
\usepackage[T1]{fontenc}

\newcommand{\compositionalityIssue}{3}
\newcommand{\compositionalityVolume}{8}
\newcommand{\paperdoinumber}{10.46298/compositionality-\compositionalityVolume-\compositionalityIssue}
\newcommand{\paperdoi}{https://doi.org/\paperdoinumber}
\newcommand{\receiveddate}{2025-04-19}
\newcommand{\accepteddate}{2026-03-03}

\usepackage{stmaryrd}
\usepackage[long, nolistings, short]{optional}

\usepackage{color}
\usepackage{amsmath, latexsym, amsthm, amsfonts,bm,amssymb}
\usepackage{natbib}
\usepackage[dvipsnames]{xcolor}
\usepackage[colorlinks,linkcolor=blue, citecolor=blue, urlcolor=blue]{hyperref}
\usepackage{commath, algorithm}

\usepackage{algpseudocode}

\algrenewcommand\algorithmiccomment[1]{\hfill\textcolor{gray}{\(\triangleright\) #1}}

\algnewcommand\Input{\item[\textbf{Input:}]}
\algnewcommand\Output{\item[\textbf{Output:}]}

\usepackage{booktabs}

\newcommand{\forw}[1]{\cF_{ #1}}

\newcommand{\backw}[1]{\cB_{#1}}
\DeclareMathSymbol{\shortminus}{\mathbin}{AMSa}{"39}

\newcommand{\revision}[1]{#1}

\newcommand{\Stoch}{\mathrm{Stoch}}

\newcommand{\Optic}{\mathrm{Optic}}
\newcommand{\Set}{\mathrm{Set}}

\usepackage{tikz-cd}
\usetikzlibrary{positioning}
\usetikzlibrary{fit}
\usetikzlibrary{cd}
\usetikzlibrary{arrows}
\usetikzlibrary{calc}
\usetikzlibrary{decorations.markings}
\usetikzlibrary{shapes.geometric}
\tikzset{ed/.style={auto,inner sep=2pt,font=\scriptsize}}
\tikzset{>=stealth}

\tikzset{vert/.style={draw,circle, minimum size=6mm, inner sep=0pt, fill=white}}
\tikzset{vertblank/.style={ minimum size=6mm, inner sep=0pt, fill=white}}

\tikzset{vertbig/.style={draw,circle, minimum size=8mm, inner sep=0pt, fill=white}}
\tikzset{->-/.style={decoration={
      markings,
      mark=at position #1 with {\arrow{>}}},postaction={decorate}}}

\tikzset{edge/.style={line width=0.5pt, decoration={markings,mark=at position 1 with %
    {\arrow[scale=1.5,>=stealth]{>}}},postaction={decorate}}}

\tikzset{dotted/.style={black!30, line width=0.5pt}}

\pgfdeclarelayer{edgelayer}
\pgfdeclarelayer{nodelayer}
\pgfsetlayers{edgelayer,nodelayer,main}

\tikzstyle{morphism}=[fill=white, draw=black, shape=rectangle]
\tikzstyle{medium box}=[fill=white, draw=black, shape=rectangle, minimum width=0.8cm, minimum height=0.9cm]
\tikzstyle{large morphism}=[fill=white, draw=black, shape=rectangle, minimum width=1.7cm, minimum height=1cm]
\tikzstyle{bn}=[fill=black, draw=black, shape=circle, inner sep=1.5pt]
\tikzstyle{effect}=[fill=white, draw=black, regular polygon, regular polygon sides=3, minimum width=0.4cm, inner sep=0pt]
\tikzstyle{state}=[fill=white, draw=black, regular polygon, regular polygon sides=3, minimum width=0.4cm, shape border rotate=180, inner sep=0pt]
\tikzstyle{medium state}=[fill=white, draw=black, regular polygon, regular polygon sides=3, minimum width=1.3cm, inner sep=0pt, shape border rotate=180]
\tikzstyle{large state}=[fill=white, draw=black, regular polygon, regular polygon sides=3, minimum width=2.2cm, shape border rotate=180, inner sep=0pt]
\tikzstyle{wn}=[fill=white, draw=black, shape=circle, inner sep=1.5pt]

\tikzstyle{arrow}=[->]
\tikzstyle{dashed box}=[-, dashed]

\tikzset{none/.style={%
     append after command={%
       \pgfextra{\node [right] at (\tikzlastnode.mid east) {{\tiny\tikzlastnode}};}
     }}}
\tikzstyle{none}=[]

\DeclareSymbolFont{bbsymbol}{U}{bbold}{m}{n}
\DeclareMathSymbol{\bbsemi}{\mathbin}{bbsymbol}{"3B}
\DeclareMathSymbol{\bbcomma}{\mathbin}{bbsymbol}{"2C}
\newcommand{\comp}{}

\newcommand{\Dup}{{\lhd}}

\newcommand{\dd}{{\,\mathrm d}}

\renewcommand{\phi}{\varphi}
\newcommand{\scr}[1]{{\mathcal #1}}

\newcommand{\EE}{\mathbb{E}}

\newcommand{\Bm}{\begin{bmatrix}}
\newcommand{\Em}{\end{bmatrix}}

\newcommand{\id}{\operatorname{id}}

\renewcommand{\P}{\mathbb{P}}

\theoremstyle{definition}
\newtheorem{thm}{Theorem}
\newtheorem{prop}[thm]{Proposition}
\newtheorem{lem}[thm]{Lemma}

\newtheorem{rem}[thm]{Remark}
\newtheorem{ex}[thm]{Example}
\newtheorem{defn}{Definition}
\newcommand{\E}{\mathbb{E}}

\newcommand{\calM}{\mathcal{M}_{\mathrm s}}
\newcommand{\scrO}{\scr{O}}

\makeatletter

\def\lstAZ{A, B, C, D, E, F, G, H, I, J, K, L, M, N, O, P, Q, R, S, T, U, V, W, X, Y, Z}
\def\lstaz{a, b, c, d, e, f, g, h, i, j, k, l, m, n, o, p, q, r, s, t, u, v, w, x, y, z}

\def\lstAZBB{B, C, D, E, F, G, H, I, J, K, L, M, N, O, P, Q, R, T, U, V, W, X, Y, Z}

\newcommand{\MkScr}[1]{\expandafter\def\csname s#1\endcsname{\mathscr{#1}}}
\newcommand{\MkUp}[1]{\expandafter\def\csname u#1\endcsname{\mathrm{#1}}}
\newcommand{\MkBold}[1]{\expandafter\def\csname b#1\endcsname{\mathbf{#1}}}
\newcommand{\MkFrak}[1]{\expandafter\def\csname f#1\endcsname{\mathfrak{#1}}}
\newcommand{\MkCal}[1]{\expandafter\def\csname c#1\endcsname{\mathcal{#1}}}
\newcommand{\MkBB}[1]{\expandafter\def\csname #1#1\endcsname{\mathbb{#1}}}

\@for\i:=\lstAZ\do{%
	\expandafter\MkScr \i  %
	\expandafter\MkFrak \i  %
	\expandafter\MkUp \i %
	\expandafter\MkBold \i %
	\expandafter\MkCal \i  %
}    
\@for\i:=\lstaz\do{%
	\expandafter\MkUp \i   }    

\@for\i:=\lstAZBB\do{%
	\expandafter\MkBB \i     }

\makeatother

\DeclareFontFamily{U}{matha}{\hyphenchar\font45}
\DeclareFontShape{U}{matha}{m}{n}{
      <5> <6> <7> <8> <9> <10> gen * matha
      <10.95> matha10 <12> <14.4> <17.28> <20.74> <24.88> matha12
      }{}
\DeclareSymbolFont{matha}{U}{matha}{m}{n}
\DeclareMathSymbol{\varrightharpoonup}{3}{matha}{"E1}
\DeclareMathSymbol{\varleftharpoonup}{3}{matha}{"E0}

\renewcommand{\tilde}{\widetilde}

\newsavebox\bsbcopier
\savebox\bsbcopier{%
  \begin{tikzpicture}[baseline=0pt,line width=0.5pt]
    \node[bn,scale=0.7] (a) at (0, 2.8pt) {};
    \draw (a) -- +(-180:.21);
    \draw (a) -- +(-45:.21);
    \draw (a) -- +(45:.21);
  \end{tikzpicture}} 
\newsavebox\bsbcopierb
\savebox\bsbcopierb{%
  \begin{tikzpicture}[baseline=0pt,line width=0.5pt]
    \node[bn,scale=0.7] (a) at (0, 2.8pt) {};
    \draw (a) -- +(0:.21);
    \draw (a) -- +(-225:.21);
    \draw (a) -- +(-135:.21);
  \end{tikzpicture}}

\newsavebox\bsinver
\savebox\bsinver{%
  \begin{tikzpicture}[baseline=0pt,line width=0.5pt]
    \node[wn,scale=0.7] (a) at (0, 2.8pt) {};
    \draw (a) -- +(0:.21);
    \draw (a) -- +(-180:.21);
  \end{tikzpicture}}

\newcommand{\bw}[0]{\mathrm{sm}}
\newcommand{\fw}[0]{\mathrm{fm}}

\usepackage{mathtools}

\DeclarePairedDelimiterXPP\expec[1]{\E}[]{}{

#1}

\begin{document}

\title{Compositionality in algorithms for smoothing}
\date{}

\author{Moritz Schauer}
\email{smoritz@chalmers.se}
\homepage{https://math.chalmers.se/~smoritz/}
\orcid{0000-0003-3310-7915}
\affiliation{Chalmers University of Technology and University of Gothenburg, Sweden}


\author{Frank van der Meulen}
\email{f.h.van.der.meulen@vu.nl}
\homepage{https://fmeulen.github.io/}
\orcid{0000-0001-7246-8612}
\affiliation{Department of Mathematics, Vrije Universiteit Amsterdam, The Netherlands}

\author{Andi Q. Wang}
\email{andi.wang@warwick.ac.uk}
\homepage{https://warwick.ac.uk/fac/sci/statistics/staff/academic-research/wang/}
\orcid{0000-0001-9551-6592}
\affiliation{Department of Statistics, University of Warwick, UK}



\begin{abstract}
Backward Filtering Forward Guiding (BFFG) is a bidirectional algorithm used for Bayesian inference on partially observed systems, first proposed in \cite{mider2021continuous} and further studied in \cite{van2020automatic}. In category theory, optics have been proposed for modelling systems with bidirectional data flow. We connect BFFG with optics by demonstrating that the forward and backwards map together define a functor from a category of Markov kernels into a category of optics, which is furthermore lax monoidal in the case when the guiding functions and kernels used in the backward step coincide with the generative dynamics. 
\end{abstract}

\maketitle

\numberwithin{equation}{section}
\sloppy

\section{Introduction}

With recent advances in computing power and hardware architectures, statistical and machine-learning models and algorithms have increased in scale and complexity. A key technique to understand and analyse complex algorithms is to understand their compositional structure: what are the fundamental building blocks that constitute an algorithm, and how are they related?

Mathematically, the study of structure and composition is the essence of \textit{category theory} \cite{pierce1991basic}. As such, the past decade has witnessed an increased application of category theory to machine learning, e.g.~\cite{Fritz2024}, \cite{shiebler2021category}, \cite{fong2019backprop}. 
One key underlying idea is to focus on compositionality and abstract away specific implementational details. 

In this work, we will make use of tools from category theory in order to study an algorithm for Bayesian inference on partially observed systems.    
Such systems include \textit{state space models}, also known as hidden Markov models, which consist of a latent Markov process that is only partially observed, possibly subject to noise. A visualisation as a directed graph is given in Figure \ref{fig:ssm}; $\bullet$ and $\circ$ correspond to latent and observed vertices respectively. 
With $X_s$ denoting the value at vertex $s$, the values at the leaf vertices are generated as follows:
\begin{align*}
	X_{t_0} \mid X_{r} &\sim \kappa_{r,t_0}(X_r,\cdot) ;\\  X_{t_i} \mid X_{t_{i-1}} & \sim \kappa_{t_{i-1}, t_i}(X_{t_{i-1}},\cdot), \qquad 1\le i \le n;  \\  X_{v_i} \mid X_{t_i} & \sim \kappa_{v_i, t_i}(X_{t_i},\cdot).
\end{align*}
Here, the $(\kappa_{s, t})$ are Markov kernels (Cf.\ Section \ref{sec:mk}) that represent the conditional probability distributions by which the process evolves on the directed graph. 
State-space models have widespread use in many fields including engineering, economics and biology (see e.g.\ \cite{rabiner2007theory}, \cite{west1997bayesian} and \cite{durbin2012time}).

\begin{figure}
\begin{center}
\begin{tikzpicture}[style={scale=0.52}]
	\tikzstyle{empty}=[fill=white, draw=black, shape=circle,inner sep=1pt, line width=0.7pt]
	\tikzstyle{solid}=[fill=black, draw=black, shape=circle,inner sep=1pt,line width=0.7pt]
	\begin{pgfonlayer}{nodelayer}
		\node [style=empty,label=below:{$r$},] (00-) at (-8, 0) {};
		\node [style=solid,label=below:{$t_0$},] (00) at (-6, 0) {};
		\node [style=empty,label={$v_0$},] (0obs) at (-6, 1.5) {};

		\node [style=solid,label=below:{$t_1$},] (0) at (-4, 0) {};
		\node [style=empty,label={$v_1$},] (1obs) at (-4, 1.5) {};

		\node [style=solid,label=below:{$t_2$}] (1) at (-2, 0) {};
		\node [style=empty,label={$v_2$},] (2obs) at (-2, 1.5) {};

		\node [] (2) at (-0, 0) {};

		\node [style=none] (end) at (1.0, 0) {.};

		\node [style=solid,label=below:{$t_{n-1}$}] (nn) at (2, 0) {};
		\node [style=empty,label={$v_{n-1}$},] (nnobs) at (2, 1.5) {};

		\node [style=solid,label=below:{$t_n$}] (n) at (4, 0) {};
		\node [style=empty,label={$v_{n}$},] (nobs) at (4, 1.5) {};

		\node [style=none] (end) at (1.0, 0) {};

	\end{pgfonlayer}
	\begin{pgfonlayer}{edgelayer}
		\draw [style=edge,color=red] (00-) to (00);
		\draw [style=edge] (00) to (0);
		\draw [style=edge] (0) to (1);
		\draw [style=edge] (1) to (2);
		\draw [style=edge] (nn) to (n);

		\draw [style=edge,color=blue] (00) to (0obs);
		\draw [style=edge,color=blue] (0) to (1obs);
		\draw [style=edge,color=blue] (1) to (2obs);

		\draw [style=edge,color=blue] (nn) to (nnobs);
		\draw [style=edge,color=blue] (n) to (nobs);

		\draw [style=dashed box] (2) to (nn);
	\end{pgfonlayer}
\end{tikzpicture}

\caption{Visualisation of a state-space-model. The process starts from the known root vertex $r$ and evolves over times $t_0, t_1,\ldots, t_n$. At each time $t_i$, a partial observation at vertex $v_i$ is assumed. \label{fig:ssm}}
\end{center}
\end{figure}

\begin{figure}
\begin{center}
\begin{tikzpicture}

\tikzstyle{empty}=[fill=white, draw=black, shape=circle,inner sep=1pt, line width=0.7pt]
\tikzstyle{solid}=[fill=black, draw=black, shape=circle,inner sep=1pt,line width=0.7pt]

\begin{pgfonlayer}{nodelayer}
		\node [style=empty,label={$r$},] (r) at (-8.5, 2) {};
		\node [style=solid,label={$t_1$},] (t1) at (-6.5, 2) {};
		\node [style=solid,label={$t_2$},] (t2) at (-4.5, 2) {};
		\node [style=solid,label={$t_4$},] (t4) at (-3, 1.5) {};
		\node [style=solid,label={$t_3$},] (t3) at (-3, 2.5) {};
		\node [style=empty,label={{ $v_1$}},] (v1) at (-1, 2.5) {};
		\node [style=empty,label={{ $v_2$}},] (v2) at (-1, 1.5) {};

\end{pgfonlayer}
\begin{pgfonlayer}{edgelayer}
		\draw [style=edge, color = red ] (r) to (t1);
		\draw [style=edge] (t1) to (t2);
		\draw [style=edge] (t2) to (t3);	
		\draw [style=edge] (t2) to (t4);
		\draw [style=edge, color=blue] (t3) to (v1);
		\draw [style=edge, color=blue] (t4) to (v2);
\end{pgfonlayer}

\end{tikzpicture}
\caption{Example of a stochastic process on a directed tree. The root-node is denoted by $r$ which, without loss of generality, is assumed to be known. We assume at each vertex sits a random quantity, and only at the leaf vertices $v_1$ and $v_2$ a realisation of that random quantity is observed.\label{fig:directed_tree}}
\end{center}
\end{figure}

In this paper we consider the more general setup where  $X$ is a stochastic process on a directed tree. The process is assumed to  be observed exactly at its leaf vertices. A simple example  is given in Figure \ref{fig:directed_tree}.

The smoothing problem  consists of inferring the distribution of the process at all interior (non-observed) vertices, conditional on the observations at the leaves. For the example in Figure \ref{fig:directed_tree} this amounts to $\{X_{t_i},\, 1\le i \le 4\}$. \cite{cappe2005springer} present multiple approaches towards this problem for the specific setting of state-space models: {\it (i)} normalised forward-backward recursion;
{\it (ii)} forward decomposition; {\it (iii)} backward decomposition. The computations involved can only be performed in closed-form in very specific settings, such as when the latent states take values only in a finite set. The celebrated Kalman filter-smoother for linear Gaussian systems is another well-known example (\cite{Bishop07}).  

\cite{mider2021continuous} introduced the {\it Backward Filtering Forward Guiding (BFFG)} algorithm (Algorithm~\ref{alg:bffg}), which uses the forward decomposition with  simplified dynamics in the backward filtering pass for $X$.   BFFG provides a general framework for dealing with the smoothing problem, encompassing both the Kalman and finite-state space settings as special cases. 
When no separate simplifying dynamics are used, the algorithm reduces to \textit{Backward Filtering Forward Sampling} (BFFS), Algorithm~\ref{alg:bffs}. The well-known \textit{Forward Filtering Backward Sampling} (FFBS) algorithm for inference in state space models (\cite{fruhwirth1994data}, \cite{carter1994gibbs}) can be seen as BFFS applied to the \textit{time-reversed} process.   

The BFFG algorithm is characterised by a {\it bidirectional data flow}: there is both a forward and backward pass. Within category theory, \textit{lenses} and \textit{optics} have been introduced as abstract models for bidirectionality; see \cite{Riley2018}.  
It is the aim of this paper to provide a categorical perspective on the BFFG algorithm using optics.
In particular, in Theorems~\ref{thm:functor} and \ref{thm:lax} we provide compositionality results for smoothing. More precisely, we show that the BFFG algorithm defines a {functor} from a category of Markov kernels into a category of optics, which in the case of BFFS is \textit{lax monoidal}. In other words, the optic obtained from parallel or sequential composition of Markov kernels is equivalent to the corresponding composition of optics of the separate Markov kernels. 
Since FFBS (\cite{fruhwirth1994data}) can be seen as a special case of BFFS, analogous compositionality results for FFBS are a direct consequence of our results.
{While BFFG in general fails to be lax monoidal, this does not typically affect its practical application in Monte Carlo inference, see Remark~\ref{rem:failure lax}.

This categorical perspective enables prospective users to break down this complex algorithm into its constituent parts, to aid implementation and comprehensibility. It demonstrates that there is a higher-order structure-preserving core within such smoothing algorithms. 
In particular, our results show that it suffices to define forward and backward maps over each edge in the graph (as in Definitions~\ref{def:backwardmap} and \ref{def:forwardmap} in Section~\ref{sec:cat_bffg}). Upon pairing both maps on an edge via the optic, the smoothing algorithm on a directed tree then consists of a composition of optics. 
Several illustrations of this are discussed in Section~\ref{subsec:examples}.

\subsection{Related work}
There have been a number of recent papers seeking to utilise concepts from category theory to study, understand and develop algorithms arising in computational statistics and machine learning.
To this end, the development of \textit{categorical probability} has been very fruitful; see for instance, \cite{Fritz2020, Cho_2019}. This approach sees Markov kernels and their composability via Chapman--Kolmogorov as primary, and enables the use of powerful tools from category theory such as string diagrams. In this work, we will use the fact that Markov kernels form a category, recalled in Proposition~\ref{prop:stoch_cat}.

A notable recent work is \cite{Fritz2024} (see also the parallel work of \cite{Virgo2023}), which uses the framework of Markov categories to abstract and generalize methods for inference on hidden Markov models (HMMs). This enables the authors to provide a unified treatment of several such smoothing algorithms. In contrast with the present work, we are less interested in studying algorithms for HMMs in the abstract, but focus on the Backward Filtering Forward Guiding algorithm in particular -- which can be applied more generally than the HMM setting -- and in fact show that this algorithm itself possesses a compositional structure.

Optics as models of bidirectional data flow have also found numerous applied uses. For instance, in the recent works of \cite{Smithe2020, Smithe2021, Braithwaite2023}, a category of optics termed \textit{Bayesian lenses} is used to give a compositional account of Bayesian inference itself. Further detailed comparisons with these particular works are given in Section~\ref{sec:cat_bffg}.
Optics have also been applied in various other places \cite{gavranovic2022space}, such as the Haskell library \textit{lens} and in the study of reverse-mode automatic differentiation \cite{Vertechi_dependent_optics}.

\subsection{Outline}
In Section \ref{sec:bffg_intro} we recap some results on Markov kernels and present the BFFG algorithm. In Section \ref{sec:cat} we recall the necessary definitions from category theory in order to state our results; specific attention is paid to the category of optics. The categorical results on BFFG together with the main results of the paper, Theorem \ref{thm:functor} and Theorem \ref{thm:lax}, are in Section \ref{sec:cat_bffg}. Section~\ref{subsec:examples} concludes the paper with several examples.

\bigskip

{\bf Acknowledgment:} We thank Keno Fischer,  Philipp Gabler,  Evan Patterson, Chad Scherrer and Bruno Gavranovi\'c for helpful discussions. AQW was partially supported by a Discipline Hopping Award from the Prob\_AI Hub (EP/Y028783/1). We are very grateful to the associate editor and anonymous referees for their extremely helpful comments, including the suggestion of the first/second terminology.

\section{Backward Filtering Forward Guiding}\label{sec:bffg_intro}
 In this section, we describe the BFFG algorithm introduced in \cite{mider2021continuous} and  generalised in \cite{van2020automatic}; see also \cite{bffg_introduction} for a non-technical introduction. We first recap standard results on Markov kernels. 
 
\subsection{Markov kernels}\label{sec:mk}
Let $S=(E,\fB)$ and $S'=(E',\fB')$ be 
 measurable spaces. A Markov kernel from source $S$ to target $S'$ is denoted by $\kappa\colon S \rightarrowtriangle S'$  (note the specific arrow). That is, $\kappa\colon E \times \fB' \to [0,1]$, where {\it (i)} for fixed $B\in \fB'$ the map $x\mapsto \kappa(x, B)$ is $\fB$-measurable and {\it (ii)} for fixed $x\in E$, the map $B\mapsto \kappa(x, B)$ is a probability measure on $S'$. 

On a measurable space $S$, denote by $\mathcal M(S)$ and $\mathbf{B}(S)$ the sets of nonnegative measures and bounded measurable functions on $S$  respectively. 
The kernel $\kappa$ induces a {\it pushforward} of a measure $\mu$ on $S$ to a measure on $S'$   via 
\begin{equation}\label{eq:pushforward}
\mu \kappa(\cdot)  = \int_{E'} \kappa(x, \cdot) \mu(\!\dd x), \qquad \mu \in \mathcal M(S).
\end{equation}
We can define a linear operator, continuous with respect to the topologies generated by the supremum norm, $\kappa\colon \mathbf B(S') \rightarrow \mathbf B(S)$ given by
\begin{equation}\label{eq:pullback}
\kappa h(\cdot) = \int_{E} h(y) \kappa(\cdot, \dd y),\qquad  h \in \mathbf{B}(S').
\end{equation} 
We will refer to this operation as the {\it pullback} of $h$ under the kernel $\kappa$. 
Markov kernels $\kappa_1\colon S_0 \rightarrowtriangle S_1 $ and $\kappa_2\colon S_1 \rightarrowtriangle S_2$  
can be composed  by the Chapman--Kolmogorov equations to obtain a kernel $S_0\rightarrowtriangle S_2$, this (associative) composition is here written as the juxtaposition
\begin{equation}\label{eq:chapman}
\kappa_2\circ \kappa_1 (x_0,\cdot) =(\kappa_1 \comp \kappa_2)(x_0, \cdot) = \int_{E_1}  \kappa_2(x_1, \cdot) \kappa_1(x_0, \dd x_1),\qquad x_0 \in E_0.
\end{equation}
We can furthermore define for $\kappa\colon S_1 \rightarrowtriangle S_2$, $\kappa'\colon S'_1 \rightarrowtriangle S_2'$
the product kernel $\kappa \otimes \kappa'\colon S_1\otimes S_1'\rightarrowtriangle S_2\otimes S_2'$ 
on the cylinder sets of $S_2 \otimes S_2'$ by 
\begin{equation}\label{eq:markovpar}
(\kappa  \otimes \kappa') ((x, x'), B \times B') = \kappa (x, B) \kappa' ( x', B'), \end{equation}  	
where $x\in E_1$, $x' \in E_1'$, $B \in \fB_2$, $B' \in \fB'_2$,
and then extend it to a kernel on the product measure space.

\subsection{Description of the stochastic dynamics on a tree}
Using Markov kernels, we can define a ``computational graph'' of the stochastic process evolving on a  \textit{directed tree}. This is most easily done through an example. Consider the  tree depicted in Figure \ref{fig:directed_tree}. Though this visualisation is popular in the literature on Bayesian networks -- highlighting vertices at which random quantities ``sit''-- from a categorical point of view it is more natural to convert this figure into the \textit{string diagram} given in Figure \ref{fig:string_directed_tree}; {for other papers also representing Bayesian networks using string diagrams, see \cite{Fong2013, Jacobs2019,Fritz2023dsep,  Fritz2023}.} {The string diagram corresponding to a state-space model is shown in Figure \ref{fig:string_ssm}, where for simplicity we took $n=2$.}

Here, for any measurable space $S$, the copy kernel $\Dup\colon S\rightarrowtriangle S\otimes S$ is a Markov kernel defined by
\begin{equation}\label{eq:dup}
\Dup(x, \dd y) = \delta_{x}(\dif y_1) \delta_x(\dif y_2), \qquad y = (y_1, y_2). 	
\end{equation}
This Markov kernel encodes the process of making a copy of $x$ when the process branches into conditionally independent sub-branches. 
The ``computational graph'' of the dynamics along the tree is then given by 
\[  \kappa_{r,t_1} \kappa_{t_1,t_2}\Dup\, (\kappa_{t_2,t_3} \otimes \kappa_{t_2,t_4}) (\kappa_{t_3,v_1} \otimes \kappa_{t_4,v_2}) . \]
By considering such generations, we see that the sequential structure of this computational graph can be simplified to
$\kappa_1 \kappa_2 \kappa_3 \kappa_4 \kappa_5$, 
where $\kappa_1=\kappa_{\tau, t_1}$, $\kappa_2=\kappa_{t_1,t_2}$, $\kappa_3=\Dup$, $\kappa_4=\kappa_{t_2,t_3} \otimes \kappa_{t_2,t_4}$ and $\kappa_5=\kappa_{t_3,v_1} \otimes \kappa_{t_4,v_2}$. 
More generally, for a directed tree, we sequentially order the vertices, by partitioning the vertices in \emph{generations} $\Gamma_1, \Gamma_2, \dots, \Gamma_n $ and $\Gamma_0 = \{0\}$, $\Gamma_{n+1} = \cV$, such that all parents of vertices $s \in \Gamma_i$ are in $\bigcup_{0 \le j < i} \Gamma_j$. 
Write $X_i = X_{\Gamma_i}$ and set
\[
\kappa_i(x_{i-1}; \dd x_{i}) = \P(X_i \in \dd x_i \mid X_{i-1} = x_{i-1}).
\]

\begin{figure}
\begin{center}
\begin{tikzpicture}

\tikzstyle{empty}=[fill=white, draw=black, shape=circle,inner sep=1pt, line width=0.7pt]
\tikzstyle{solid}=[fill=black, draw=black, shape=circle,inner sep=1pt,line width=0.7pt]

\begin{pgfonlayer}{nodelayer}
		\node [style=empty,] (r) at (-10, 2) {};
		\node [style=solid,] (t1) at (-8, 2) {};
		\node [style=morphism, draw=black] (t2) at (-4.5, 2) {$\Dup$};
		\node [style=solid,] (t4) at (-2, 1) {};
		\node [style=solid,] (t3) at (-2, 3) {};
		\node [style=empty,] (v1) at (0, 3) {};
		\node [style=empty,] (v2) at (0, 1) {};
		\node [style=morphism, draw=black] (k6) at (-9, 2) {$\kappa_{r,t_1}$};
		\node [style=morphism, draw=black] (k6) at (-6.8, 2) {$\kappa_{t_1,t_2}$};
		\node [style=solid,] (k6) at (-5.75, 2) {};

		\node [style=morphism, draw=black] (k6) at (-3.25, 2.5) {$\kappa_{t_2,t_3}$};
		\node [style=morphism, draw=black] (k6) at (-1, 3) {$\kappa_{t_3,v_1}$};
		\node [style=morphism, draw=black] (k6) at (-3.25, 1.5) {$\kappa_{t_2,t_4}$};
		\node [style=morphism, draw=black] (k6) at (-1, 1) {$\kappa_{t_4,v_2}$};

\end{pgfonlayer}
\begin{pgfonlayer}{edgelayer}
		\draw [style=none, color = red ] (r) to (t1);
		\draw [style=none] (t1) to (t2);
		\draw [style=none] (t2) to (t3);	
		\draw [style=none] (t2) to (t4);
		\draw [style=none, color=blue] (t3) to (v1);
		\draw [style=none, color=blue] (t4) to (v2);

\end{pgfonlayer}

\end{tikzpicture}
\caption{String diagram corresponding to Figure \ref{fig:directed_tree}. $\Dup$ is is the duplication kernel defined in \eqref{eq:dup}. \label{fig:string_directed_tree}}
\end{center}
\end{figure}

\begin{figure}
\begin{center}
\begin{tikzpicture}[style={scale=0.64}]
    \tikzstyle{empty}=[fill=white, draw=black, shape=circle,inner sep=1pt, line width=0.7pt]
    \tikzstyle{solid}=[fill=black, draw=black, shape=circle,inner sep=1pt,line width=0.7pt]
	\begin{pgfonlayer}{nodelayer}
        \node [style=empty] (r)    at (-3.0, 1) {};
        \node [style=morphism, draw=black] (k1)   at (-1.5, 1)  {$\kappa_{r,t_0}$};
        \node [style=solid] (s1)   at (0.0, 1)   {};
        \node [style=morphism, draw=black] (dup1) at (1.25, 1)  {$\Dup$};
        \node [style=morphism, draw=black] (kp1)  at (2.75, 2)  {$\kappa_{t_0,v_0}$};
        \node [style=empty] (v1)   at (4.25, 2)  {};
        \node [style=morphism, draw=black] (k2)   at (2.75, 0)  {$\kappa_{t_0,t_1}$};
        \node [style=solid] (s2)   at (4.25, 0)  {};
        \node [style=morphism, draw=black] (dup2) at (5.5, 0)   {$\Dup$};
        \node [style=morphism, draw=black] (kp2)  at (7.0, 1)   {$\kappa_{t_1,v_1}$};
        \node [style=empty] (v2)   at (8.5, 1)   {};
        \node [style=morphism, draw=black] (k3)   at (7.0, -1)  {$\kappa_{t_1,t_2}$};
        \node [style=solid] (s3)   at (8.5, -1)  {};
        \node [style=morphism, draw=black] (kp3)  at (10.0, -1) {$\kappa_{t_2,v_2}$};
        \node [style=empty] (v3)   at (11.5, -1) {};
        \node [style=none] (end)   at (10.0, -2.5) {};
        \node [style=none] (endr)  at (11.5, -2.5) {};
	\end{pgfonlayer}
	\begin{pgfonlayer}{edgelayer}
        \draw [style=none, color=red] (r) to (k1);
        \draw [style=none, color=red] (k1) to (s1);
        \draw [style=none] (s1) to (dup1);
        \draw [style=none, color=blue, bend left=30] (dup1) to (kp1);
        \draw [style=none, color=blue] (kp1) to (v1);
        \draw [style=none, bend right=30] (dup1) to (k2);
        \draw [style=none] (k2) to (s2);
        \draw [style=none] (s2) to (dup2);
        \draw [style=none, color=blue, bend left=30] (dup2) to (kp2);
        \draw [style=none, color=blue] (kp2) to (v2);
        \draw [style=none, bend right=30] (dup2) to (k3);
        \draw [style=none] (k3) to (s3);
        \draw [style=none, color=blue] (s3) to (kp3);
        \draw [style=none, color=blue] (kp3) to (v3);
	\end{pgfonlayer}
\end{tikzpicture}\\
\end{center}
\caption{{String diagram corresponding to Figure \ref{fig:ssm} when $n=2$.  \label{fig:string_ssm}}}
\end{figure}

Therefore, in the following, we assume a stochastic process on a line-graph, with vertices denoted by $\{0,1,\ldots, n+1\}$, where the root-vertex $0$ is assumed to be known and the transition from vertex $i-1$ to $i$ is governed by the Markov kernel $\kappa_i:E_{i-1}\to E_i$. Then, if at node $0$ the process equals the known quantity $\hat x_0$, the distribution at vertex 1 is $\delta_{\hat x_0} \kappa_1$, the pushforward of $\delta_{\hat x_0}$ (cf.\ Equation \eqref{eq:pushforward}). Pushing forward iteratively, we see that 
the random quantity at vertex $i$, denoted by $X_i$, has  distribution
\begin{equation}\label{eq:linegraph} \delta_{\hat x_0} \kappa_1 \kappa_2 \cdots \kappa_i. \end{equation}

\subsection{Conditioning on a line graph}
Suppose that we observe $x_{n+1}$ at vertex $n+1$.
See Figure \ref{fig:linegraph}, where we have added the Markov kernels along the branches. The \textit{smoothing distribution} is the distribution of the process on the graph, conditioned on its value $x_{n+1}$ at node $n+1$.
We are interested in generating samples from the smoothing distribution.

\begin{figure}
\begin{center}
\begin{tikzpicture}[style={scale=0.52}]

\tikzstyle{empty}=[fill=white, draw=black, shape=circle,inner sep=1pt, line width=0.7pt]
\tikzstyle{solid}=[fill=black, draw=black, shape=circle,inner sep=1pt,line width=0.7pt]

\begin{pgfonlayer}{nodelayer}
		\node [style=empty,label={$0$},] (n0) at (0, 0) {};
		\node [style=solid,label={$1$},] (n1) at (4, 0) {};
		\node [style=solid,label={$2$},] (n2) at (8, 0) {};
		\node [style=solid,] (intermediate) at (11, 0) {};
		\node [style=solid,label={$n$},] (n) at (14, 0) {};
		\node [style=empty,label={$n+1$},] (n+1) at (18.5, 0) {};
		\node [style=morphism] (kap1) at (2, 0){$\kappa_{1}$};
		\node [style=morphism] (kap2) at (6, 0){$\kappa_{2}$};
		\node [style=morphism] (kap2) at (16, 0){$\kappa_{n+1}$};
\end{pgfonlayer}

\begin{pgfonlayer}{edgelayer}
		\draw [style=edge] (n0) to (n1);
		\draw [style=edge] (n1) to (n2);
		\draw [style=edge] (n2) to (intermediate);
		\draw [style=edge] (n) to (n+1);
		\draw [style=dashed box] (intermediate) to (n);
\end{pgfonlayer}

\end{tikzpicture}

\caption{Stochastic process on a line graph with one observation corresponding to the composition of Markov kernels specified in Equation \eqref{eq:linegraph}. \label{fig:linegraph}}
\end{center}
\end{figure}

We assume that for the final vertex {the kernel $\kappa_{n+1}(x,\dif y)$ has a jointly measurable density {$p_{n+1}(x,y)$} with respect to a dominating measure $\nu$:
\begin{equation*}
    \kappa_{n+1}(x,A)=\int_A p_{n+1}(x,y)\nu(\dif y).
\end{equation*}
} 
We define the family of functions $(h_i)_{i=1}^n$ as follows.
\begin{defn}\label{defn:backwardfilter}
On a line graph, the \emph{backward filter} is defined by the following recursive scheme, where we recurse from $i=n$ \textit{down} to $i=1$:
\begin{equation}\label{eq: hn2}
\begin{split}
h_n(\cdot) &= p_{n+1}(\cdot, x_{n+1}),\\
h_{i-1} & = \kappa_{i} h_{i}, \qquad i=1,\ldots, n,
\end{split}
\end{equation}
each $h_{i-1}$ defined as a pullback as in \eqref{eq:pullback}.    
\end{defn}
Each $h_i(x)$ is the \textit{likelihood} corresponding to a fixed observation
$X_{n+1}=x_{n+1}$, given `parameter' $X_i=x$. 
The stochastic process that is  conditioned on $X_{n+1} = x_{n+1}$ is Markovian then with transition kernels
\begin{equation}\label{eq:kappa_to_kappastar}
    \kappa^\star_i(x,\dd y) = \frac{\kappa_i(x,\dd y) h_i(y)}{\int_{E_i} \kappa_i(x,\dd y') h_i(y')}. 
\end{equation} 
This motivates the following definition.
\begin{defn}\label{def:htransform}
For a family of strictly positive functions $(h_i)$, the  {\it $h$-transform} maps each  $(\kappa_i, h_i)$ to $\kappa_i^\star$, where $\kappa_i^\star$ is defined in \eqref{eq:kappa_to_kappastar}.
\end{defn}

Therefore, if $(h_i)_{i=0}^n$ are known, we can pushforward {$\delta_{\revision{\hat x_0}}$, via
\begin{equation}\label{eq:linegraph*} \delta_{\hat x_0} \kappa_1^\star \kappa_2^\star \cdots \kappa_i^\star, \end{equation}
which gives the marginal distribution of the conditioned process at vertex $i$, with the joint distribution of $(X_0^\star, \dots, X_n^\star)$ being
\[ P^\star(\dif x_0,\ldots, \dd x_n) = \delta_{\hat x_0}(\dif x_0) \kappa^\star_1(x_0, \dd x_1) \cdots \kappa^\star_{n}(x_{n-1},  \dd x_{n}). \]

We remark that,
with marginal distributions of $X_i$ given by \eqref{eq:linegraph} and here denoted by $\mu_i$,
for almost all observations  (with respect to the marginal distribution of $X_{n+1}$), the functions 
$h_i(x)$  will be different from 0 on sets that have positive marginal probability (with respect to $\mu_i$). 
Therefore $\kappa^\star_i$ is well-defined up sets of zero marginal probability.

We conclude that a realisation of the conditioned process can be obtained as follows: first compute the \textit{backward} filter, to obtain the family $(h_i)_{i=0}^n$, then \textit{forward} simulate from the kernels $\kappa_i^\star$, where $\kappa_i^\star$ is the $h$-transform of $(\kappa_i, h_i)$. Hence, the operation is \textit{bidirectional}. 
This is spelt out explicitly in the Backward Filtering Forward Sampling algorithm detailed in Algorithm~\ref{alg:bffs}.

\begin{algorithm}
\caption{Backward Filtering Forward Sampling \label{alg:bffs}}
\begin{algorithmic}[1]
\Input Initial state $x_0$; observation time $n + 1$;
transition kernels $\kappa_i$;
observation $x_{n+1}$ of $X_{n+1}$; conditional $\nu$-density $p_{n+1}$.
\Output A sample path $(x_0, x_1, \dotsc, x_{n})$ from $P^\star$ and likelihood $h_0(x_0)$.

\vspace{0.5em}

\Function{BackwardFilter}{$(\kappa_1, \dots, \kappa_{n}), h_{n}, n$}
    \For{$i = n$ down to $1$}
        \State $h_{i-1} \gets \kappa_{i} h_{i} $
    \EndFor
    \State $\mathfrak m \gets (h_1, \dotsc, h_n)$
    \State \Return $h_0, \mathfrak m$
\EndFunction

\vspace{0.5em}

\Function{ForwardSample1}{$x, \kappa, h$}
    \State Let $\pi (\dif y) \propto \kappa(x, \dif y) h(y)$
    \State Sample $y \sim \pi$
    \State \Return $y$
\EndFunction
\vspace{0.5em}
\Function{ForwardSample}{$x_0, (\kappa_1,\dots, \kappa_{n}), \mathfrak m, n$}
\State $(h_1, \dots, h_n) \gets \mathfrak m$
    \For{$i = 1$ to $n$}
        \State $x_i \gets$ \Call{ForwardSample1}{$x_{i-1}, \kappa_{i}, h_{i}$}
    \EndFor
    \State \Return $(x_0, \dots, x_n)$
\EndFunction

\vspace{0.5em}

\Function{BFFS}{$x_0, (\kappa_1, \dots, \kappa_{n}), p_{n+1}, x_{n+1}, n$}
    \State $h_n(\cdot) \gets p_{n+1}(\cdot, x_{n+1})$
    \State $ h_0, \mathfrak m\gets$ \Call{BackwardFilter}{$(\kappa_1, \dots, \kappa_{n}), h_n, n$}
    \State \Return \Call{ForwardSample}{$x_0, (\kappa_1, \dots, \kappa_{n}), \mathfrak m, n$}, $h_0(x_0)$
\EndFunction
\end{algorithmic}
\end{algorithm}

\subsection{Guiding on a line graph}
In most cases of practical interest, the backward filter \eqref{eq: hn2} cannot be computed in closed form. For that reason, we associate to each kernel $\kappa_i$ an approximate   \textit{guiding} kernel $\tilde\kappa_i$ which is chosen such that the backward filter with the kernels $\tilde\kappa_i$ is tractable (that is, more easily computed, preferably in closed form): 
\begin{equation}\label{eq: g}
\begin{split}
g_n(\cdot) &= \tilde p_{n+1}(\cdot, x_{n+1})\\
g_{i-1} & = \tilde\kappa_{i} g_{i}, \qquad i=1,\ldots, n.
\end{split}
\end{equation}
{We will refer to these functions $\{g_i\}_{i=0}^n$ as the \textit{guiding functions}.}
Using our guiding functions, we define the Markov kernels obtained by the $h$-transform of $(\kappa_i, g_i)$:
\begin{equation}
    \kappa^\circ_i(x,\dd y) = \frac{\kappa_i(x,\dd y) g_i(y)}{\int_{E_i} \kappa_i(x,\dd y') g_i(y')}.
    \label{eq:kappa_circ}
\end{equation}

If chosen properly (see Remark \ref{rem:choice_g}), the sequence of maps  $\{g_i\}_{i=0}^n$ ensures that realisations of the process under $\kappa^\circ_i$ are {\it guided} towards $x_{n+1}$ at vertex $n+1$.
{Indeed, these kernels $\kappa^\circ_i$ define a \textit{guided process} with law $P^\circ$ \eqref{eq:P_circ} with the property that our observed point $x_{n+1}$ has high probability under $P^\circ$.}
{By taking quotients of \eqref{eq:kappa_to_kappastar} and \eqref{eq:kappa_circ}}, we obtain Radon--Nikodym derivatives
\[ \frac{\kappa^\star_i(x,\dd y)}{\kappa^\circ_i(x,\dd y)} = \frac{h_i(y)}{g_i(y)} \frac{\int_{E_i}\kappa_i(x,\dd y') g_i(y')}{\int_{E_i}\kappa_i(x,\dd y') h_i(y')}.\]

\begin{prop}
    The likelihood ratio of paths samples under $\kappa^\star_i$ relative to $\kappa_i^\circ$ is given by
    \begin{equation*}
        L(x_1,\ldots, x_n) =\frac{h_{n}(x_{n})}{g_n(x_n)}\frac{g_0(\hat x_0)}{h_0(\hat x_0)}\prod_{i=1}^{n} \frac{ \kappa_ig_i(x_{i-1})}{\tilde\kappa_{i}g_i (x_{i-1})}.
    \end{equation*}
    \label{prop:LR}
\end{prop}
\begin{proof}
    The likelihood ratio is
\begin{align*}
    L(x_1,\ldots, x_n) &= \prod_{i=1}^{n} \frac{h_i(x_{i})}{g_i(x_{i})} \frac{\int \kappa_i(x_{i-1},\dd y) g_i(y)}{\int \kappa_i(x_{i-1},\dd y) h_i(y)}\\ &= \frac{h_{n}(x_{n})}{h_0(\hat x_0)}\prod_{i=1}^{n} \frac{\int \kappa_i(x_{i-1},\dd y) g_i(y)}{g_i(x_i)} \\
    &= \frac{h_{n}(x_{n})}{g_n(x_n)}\frac{g_0(\hat x_0)}{h_0(\hat x_0)}\prod_{i=1}^{n} \frac{ \kappa_ig_i(x_{i-1})}{\tilde\kappa_{i}g_i (x_{i-1})} ,
\end{align*}  
where the first equality is simply the definitions, the second equality is due to the cancellation of terms due from the recursive relation in \eqref{eq: hn2} and the final equality from the recursive relation \eqref{eq: g}. 
\end{proof}

Proposition~\ref{prop:LR} implies that for bounded functions $f$,
\begin{equation}\label{eq:condexp}
    \EE \left[f(X_1,\ldots, X_n) \mid X_{n+1}=x_{n+1}\right] = \EE^\circ \left[f(X_1,\ldots, X_n) L(X_1,\ldots, X_n) \right],  
\end{equation} 
where $\EE^\circ$ denotes expectation under which $X$ has law
\begin{equation}
    P^\circ(\dif x_0,\ldots, \dd x_n) = \delta_{\hat x_0}(\dif x_0) \kappa^\circ_1(x_0, \dd x_1) \cdots \kappa^\circ_{n}(x_{n-1},  \dd x_{n}).
    \label{eq:P_circ}
\end{equation}

Assuming $h_n=p_{n+1}(\cdot, x_{n+1})$ to be tractable and $X_0=\hat x_0$ being fixed and known, the likelihood ratio does not involve $h_1,\ldots, h_{n-1}$ which may be hard or expensive to compute. It does however depend on  $h_0(\hat x_0)$, the \emph{evidence}, which can usually not be evaluated. Fortunately,  Markov Chain or Sequential Monte Carlo methods that exploit the relation in \eqref{eq:condexp} only require evaluation of $L$ up to a fixed multiplicative constant. For instance, \eqref{eq:condexp} implies
\[
\EE \left[f(X_1,\ldots, X_n) \mid X_{n+1}=x_{n+1}\right] = \frac{\EE^\circ \left[f(X_1,\ldots, X_n) L(X_1,\ldots, X_n) \right]}{\EE^\circ \left[ L(X_1,\ldots, X_n) \right]}. \] The right-hand-side can be approximated by the \textit{self-normalised importance sampling estimator},
\begin{equation}\label{eq:selfnormalised_is} T:= \frac{\sum_{i=1}^I f(X_{1,i},\ldots, X_{n,i}) L(X_{1,i},\ldots, X_{n,i})}{\sum_{i=1}^I L(X_{1,i},\ldots, X_{n,i})},
\end{equation}
where $\{(X_{1,i}, \ldots, X_{n,i})\}_{i=1}^I$ are independent samples under {${P}^\circ$}. We see that $T$ does not depend on $h_0(\hat x_0)$.

An algorithmic description of the full Backward Filtering Forward Guiding is given in Algorithm~\ref{alg:bffg}. In the case when $\tilde \kappa=\kappa$, we recover BFFS as in Algorithm~\ref{alg:bffs}.

\begin{algorithm}
\caption{Backward Filtering Forward Guiding (BFFG)\label{alg:bffg}}
\begin{algorithmic}[1]
\Input 
Initial state $x_0$; observation time $n + 1$;
transition kernels $(\kappa_i)$, $(\tilde \kappa_i)$;
observation $x_{n+1}$ of $X_{n+1}$; conditional $\nu$-densities $p_{n+1}$ and $\tilde p_{n+1}$.
\Output A sample path $(x_0, x_1, \dotsc, x_{n})$ from $P^\circ$ and a weight $\varpi = h_0(x_0) \cdot L(x_1,\dots, x_n)$.

\vspace{0.5em}

\Function{ForwardGuiding1}{$(x, \varpi), \kappa, \tilde \kappa,g$}
        \State Let $\pi (\dif y) = \kappa (x, \dif y) g(y)/(\kappa g)(x)$
        \State Sample $y \sim \pi$
        \State $w \gets  \dfrac{(\kappa g)(x)}{(\tilde \kappa g)(x)}$
        \State \Return $y, w\cdot \varpi$
\EndFunction

\vspace{0.5em}

\Function{ForwardGuiding}{$(x_0, \varpi), (\kappa_1, \dots, \kappa_{n}), (\tilde \kappa_1, \dots, \tilde \kappa_{n}),\mathfrak m, n$}
    \State $(g_1, \dots, g_n) \gets \mathfrak m$    
    \For{$i = 1$ to $n$}
        \State $x_i, \varpi \gets$ \Call{ForwardGuiding1}{$(x_{i-1}, \varpi), \kappa_{i},\tilde \kappa_{i},   g_{i}$}
    \EndFor
    \State \Return $(x_0, \dots, x_n)$, $\varpi$
\EndFunction

\vspace{0.5em}

\Function{BFFG}{$x_0, (\kappa_1, \dots, \kappa_{n}), (\tilde\kappa_1, \dots, \tilde \kappa_{n}), p_{n+1},\tilde p_{n+1}, x_{n+1}, n$}
    \State $h_n (\cdot) = p_{n+1}(\cdot, x_{n+1})$
    \State $g_n (\cdot) = \tilde p_{n+1}(\cdot, x_{n+1})$
    \State $ g_0,  \mathfrak m \gets$ \Call{BackwardFilter}{$(\tilde\kappa_1, \dots, \tilde\kappa_{n}), g_n, n$}
     \State $(x_0, \dots, x_n), \varpi \gets$ \Call{ForwardGuiding}{$(x_0, g_0(x_0)), (\kappa_1, \dots, \kappa_{n}),(\tilde \kappa_1, \dots, \tilde \kappa_{n}), \mathfrak m, n$}
    \State \Return $(x_0, \dots , x_n), \varpi\cdot h_n(x_n)/g_n(x_n)$
\EndFunction
\end{algorithmic}
\end{algorithm}

\begin{rem}\label{rem:choice_g}
We need to choose the guiding functions $
\{g_i\}$ such that  $P^\star$ is absolutely continuous with respect to $P^\circ$. A specific choice can be made for example by minimising the (reverse) Kullback--Leibler divergence between $P^\star$ and $P^\circ$; other discrepancy measures between the two measures are possible. Alternatively, if one uses BFFG within a sequential Monte Carlo algorithm, a practical criterion consists of choosing the guiding functions such that the effective sample size of particles is maximised. 
\end{rem}

\begin{rem}
For ease of exposition, we have introduced the BFFG algorithm on a line graph with a single observation. 
We refer to Theorem 2.2 in \cite{van2020automatic}  for a more general statement on a directed acyclic graph.
\end{rem}

\section{Definitions from category theory}\label{sec:cat}

The main contribution of this work is to connect the BFFG algorithm described above with well-established concepts from category theory.
Our main contribution will be showing that the BFFG algorithm \textit{preserves sequential composition} and furthermore \textit{parallel composition} in the case of BFFS.
To formalise this, we will make use of well-established categorical concepts such as \textit{monoidal categories} and \textit{lax monoidal functors}.
For readers who are less familiar with these concepts, we have included full definitions in Appendix~\ref{app:cat}.

\subsection{The category Stoch}
It is a key contribution of categorical probability, \cite{Fritz2020, Cho_2019}, that
the Chapman--Kolmogorov equations \eqref{eq:chapman}, actually constitute composition $\kappa_2 \circ \kappa_1$ in a particular category, $\Stoch$, of Markov kernels. 
Formally, we have the following well-known fact; see, for instance \cite[Section~4]{Fritz2020}:
\begin{prop}\label{prop:stoch_cat}
    The category $\Stoch$, with objects measurable spaces, $S,S',\dots$, and Markov kernels $\kappa\colon{}S\rightarrowtriangle S'$ as morphisms is a symmetric monoidal category: sequential composition is as defined in \eqref{eq:chapman}.  The identity $\id$ for this composition is the identity function considered as a Markov kernel, $\id(x, \dd y) = \delta_{x}(\!\dd y)$.
    The monoidal product is as in \eqref{eq:markovpar} with monoidal unit $1_\Stoch = \left(\{*\}, (\emptyset, \{*\})\right)$.

\end{prop}

\subsection{Optics}\label{sec:optics}
Optics were introduced in \cite{Riley2018} for modelling systems with bidirectional data flow. Here, we present the formulation as given in  \cite{gavranovic2022space}.  Note that an optic models a transformation in which a \revision{first} (forward) pass is followed by a \revision{second} (backward) pass.  

\begin{defn}\label{def:optic}
Given a symmetric monoidal base category $(\scr{C}, \otimes, 1,\alpha,\lambda,\rho)$, an {\bf optic}  from a pair of objects from $\scr C$, $(A, A')$ to a pair $(B, B')$ is an equivalence class of triplets $(M, \fw, \bw)$ where 
\begin{itemize}
	\item $M$ is an object of $\scr{C}$, the type of the \textit{internal state};
	\item $\fw \colon A \to M \otimes B$, a morphism referred to as the \textit{\revision{first} map};
	\item $\bw \colon M \otimes B' \to A'$, a morphism referred to as the \textit{\revision{second} map};
\end{itemize}
and the equivalence class is given by the finest equivalence relation containing the following relation: for two such triplets $(M, f, f')$ and $(N,g,g')$, we have $(M, f, f') \sim (N,g,g')$ if there exists a \textit{residual morphism} $r\colon M \to N$ such that 
$(r\otimes B)\circ f=g$ and  $g'\circ(r\otimes B') = f'$:
\begin{center}	
\begin{tikzcd}
A \arrow[rd, "g"] \arrow[r, "f"] & M\otimes B \arrow[d, "r\otimes B"] \\
							  & N\otimes B
\end{tikzcd}	\qquad \qquad 
\begin{tikzcd}
M\otimes B' \arrow[r, "f'"] \arrow[d, "r\otimes B'", swap] & A'  \\
							   N\otimes B' \arrow[ru, "g'", swap]
\end{tikzcd}.	
\end{center}	
\end{defn}

An optic from  $(A,A')$ to $(B, B')$ can be visualised as follows:

\begin{equation}\label{eq:optic}
	\begin{tikzpicture}[style={scale=1.5}]
	\begin{pgfonlayer}{nodelayer}
		\node [style=morphism] (bwmap) at (0, 0) {$\fw$};
		\node [style=none] (in_bw) at (-1, 0) {$A$};	
		\node [style=none] (out_bw) at (1, 0) {$B$};	
		
		\node [style=morphism] (fwmap) at (0, -1) {$\bw$};
		\node [style=none] (in_fw) at (-1, -1) {$A'$};	
		\node [style=none] (out_fw) at (1, -1) {$B'$};	

		\node [style=none] (message) at (0.25, -0.5) {$M$};	
	\end{pgfonlayer}	
	
	\begin{pgfonlayer}{edgelayer}
		\draw [style=edge] (in_bw) to (bwmap);
		\draw [style=edge] (bwmap) to (out_bw);		
		
		\draw [style=edge] (fwmap) to (in_fw);
		\draw [style=edge] (out_fw) to (fwmap);		

		\draw [style=edge] (bwmap) to (fwmap) ;		
	\end{pgfonlayer}

\end{tikzpicture}
\end{equation}
In the \revision{first} pass, $\fw$ produces a {\it message} of type $M$, the internal state, which is consumed in the \revision{second} pass by $\bw$.

\begin{rem}
    In the literature, the first pass is typically referred to as the \textit{forward} pass, and the second pass as the \textit{backward} pass, reflecting the flow of the information, e.g. \cite{gavranovic2022space}.
    We have opted to use the (new) terminology of first/second instead of forwards/backwards in order to avoid confusion: the first (i.e.\ forward) pass of the BFFG algorithm is actually the \textit{backward} filter, and the second (i.e.\ backward) pass of the BFFG algorithm is the \textit{forward} guiding step.
\end{rem}

\begin{defn}[Composition of optics]
Optics can be composed as follows \cite[Definition~3]{gavranovic2022space}: consider two optics 
\[ \begin{bmatrix}
	A \\ A'
\end{bmatrix} \xrightarrow{(M_1, \fw_1, \bw_1)}\begin{bmatrix}
	B \\ B'
\end{bmatrix} \quad \text{and} \quad  \begin{bmatrix}
	B \\ B'
\end{bmatrix} \xrightarrow{(M_2, \fw_2, \bw_2)}\begin{bmatrix}
	C \\ C'
\end{bmatrix}, \]
then the composite $(M, \fw, \bw)\colon{}(A,A')\to (C,C')$ is defined by
\begin{align}
	M & := M_1 \otimes M_2 \\
	\fw &:=   (\id_{M_1} \otimes\, \fw_2)\circ \fw_1 \colon \quad A \rightarrow M_1 \otimes (M_2 \otimes C)\label{eq:fw_comp} \\ 
	\bw &:= \bw_1 \circ (\id_{M_1} \otimes \,\bw_2) \colon \quad M_1\otimes (M_2 \otimes C') \to A',\label{eq:bw_comp}
\end{align}
\end{defn}
noting we have further isomorphisms $M_1 \otimes (M_2 \otimes C)\cong (M_1 \otimes M_2) \otimes C$ and $M_1\otimes (M_2 \otimes C')\cong (M_1\otimes M_2) \otimes C'$.
This composition is represented diagrammatically below. The blue rectangles indicate that its constituent components should be viewed together as a single unit; compare with \eqref{eq:optic}.

\begin{center}

\begin{tikzpicture}[style={scale=1.5}]
	\begin{pgfonlayer}{nodelayer}
		\node [style=morphism] (bwmap2) at (0.5, 0) {$\fw_1$};
		\node [style=none] (in_bw2) at (-1, 0) {$A$};	
		
		\node [style=morphism] (fwmap2) at (0.5, -1.5) {$\bw_1$};
		\node [style=none] (in_fw2) at (-1, -1.5) {$A'$};	
		\node [style=none] (message2) at (0.75, -0.75) {$M_1$};	
	\end{pgfonlayer}	
	\begin{pgfonlayer}{edgelayer}
		\draw [style=edge] (in_bw2) to (bwmap2);
		
		\draw [style=edge] (fwmap2) to (in_fw2);
		\draw [style=edge] (bwmap2) to (fwmap2);		
	\end{pgfonlayer}
	\begin{pgfonlayer}{nodelayer}
		\node [style=morphism] (bwmap) at (2.5, 0) {$\fw_2$};
		\node [style=none] (out_bw) at (4, 0) {$C$};	
		
		\node [style=morphism] (fwmap) at (2.5, -1.5) {$\bw_2$};
		\node [style=none] (out_fw) at (4, -1.5) {$C'$};	
		\node [style=none] (message) at (2.75, -0.75) {$M_2$};	
	\end{pgfonlayer}	
	
	\begin{pgfonlayer}{edgelayer}
		
		\draw [style=edge] (bwmap2) to (bwmap);
		
		\draw [style=edge] (fwmap) to (fwmap2);
		
		\draw [style=edge] (bwmap) to (out_bw);		
		\draw [style=edge] (out_fw) to (fwmap);		
		\draw [style=edge] (bwmap) to (fwmap);		
	\end{pgfonlayer}

	\node [style=none] (B_label) at (1.5, 0.2) {$B$};
	\node [style=none] (Bprime_label) at (1.5, -1.3) {$B'$};

    \node[draw,inner sep=2.5mm, label=below:,fit=(bwmap)  (bwmap2),color=blue] {};
    \node[draw,inner sep=2.5mm, label=below:,fit= (fwmap) (fwmap2) ,color=blue] {};
    
\end{tikzpicture}

\end{center}

\begin{defn}[Tensor product of optics]
Consider two optics 
\[ \begin{bmatrix}
	A_1 \\ A_1'
\end{bmatrix} \xrightarrow{(M_1, \fw_1, \bw_1)}\begin{bmatrix}
	B_1 \\ B_1'
\end{bmatrix} \quad \text{and} \quad  \begin{bmatrix}
	A_2 \\ A_2'
\end{bmatrix} \xrightarrow{(M_2, \fw_2, \bw_2)}\begin{bmatrix}
	B_2 \\ B_2'
\end{bmatrix}, \]
then the tensor product $(M, \fw, \bw)$, mapping
\[ \Bm A_1 \otimes A_2 \\ A_1'\otimes A_2'\Em \xrightarrow{(M, \fw, \bw)} \Bm B_1 \otimes B_2 \\ B_1' \otimes B_2'\Em \]
 is defined by
\begin{align*}
	M & := M_1 \otimes M_2 \\
	\fw & := \fw_1\otimes \,\fw_2\colon \quad A_1 \otimes A_2 \to (M_1 \otimes B_1) \otimes (M_2 \otimes B_2) \\
	\bw & :=  \bw_1 \otimes \bw_2\colon \quad (M_1\otimes B_1')\otimes (M_2\otimes B_2')   \to A_1' \otimes A_2',
\end{align*}
\end{defn}
further noting that since the base category is symmetric monoidal we have isomorphisms $(M_1 \otimes B_1) \otimes (M_2 \otimes B_2)\cong (M_1 \otimes M_2)\otimes (B_1 \otimes B_2)$ and $(M_1\otimes B_1')\otimes (M_2\otimes B_2')\cong (M_1 \otimes M_2)\otimes (B_1'\otimes B_2')$.
This tensor product can be visualised by placing the corresponding diagrams \eqref{eq:optic} adjacent and overlapping one another; again, blue rectangles indicate those sections should be viewed as a single block:

\begin{center}

\begin{tikzpicture}[style={scale=1.5}]
	\begin{pgfonlayer}{nodelayer}
		\node [style=morphism] (bwmap2) at (0.75, -0.5) {$\fw_1$};
		\node [style=none] (in_bw2) at (-1, -0.5) {$A_1$};	
		
		\node [style=morphism] (fwmap2) at (0.75, -2.5) {$\bw_1$};
		\node [style=none] (in_fw2) at (-1, -1) {$A_2$};	
		\node [style=none] (message2) at (1, -1.75) {$M_1$};
        \node[style=none] (in_Aprime1) at (-1,-2.5) {$A_1'$};
        \node[style=none] (in_Aprime2) at (-1,-3.) {$A_2'$};
	\end{pgfonlayer}	
	\begin{pgfonlayer}{edgelayer}
		\draw [style=edge] (in_bw2) to (bwmap2);
		
		\draw [style=edge] (fwmap2) to (in_Aprime1);
		\draw [style=edge] (bwmap2) to (fwmap2);		
	\end{pgfonlayer}
	\begin{pgfonlayer}{nodelayer}
		\node [style=morphism] (bwmap) at (2.25, -1) {$\fw_2$};
		\node [style=none] (out_bw) at (4, -0.5) {$B_1$};	
		
		\node [style=morphism] (fwmap) at (2.25, -3.) {$\bw_2$};
		\node [style=none] (out_fw) at (4, -1) {$B_2$};	
		\node [style=none] (message) at (2.5, -1.75) {$M_2$};
        \node [style=none] (B1prime) at (4,-2.5) {$B_1'$};
        \node [style=none] (B2prime) at (4,-3.) {$B_2'$};
	\end{pgfonlayer}	
	
	\begin{pgfonlayer}{edgelayer}
		\draw [style=edge] (bwmap2) to (out_bw);
		
		\draw [style=edge] (fwmap) to (in_Aprime2);
		
		\draw [style=edge] (bwmap) to (out_fw);		
		\draw [style=edge] (B2prime) to (fwmap);		
		\draw [style=edge] (bwmap) to (fwmap);	
        
        \draw [style=edge] (in_fw2) to (bwmap);
        \draw [style=edge] (B1prime) to (fwmap2);
	\end{pgfonlayer}

\node[draw,inner sep=2.5mm, label=below:,fit=(bwmap)   (bwmap2),color=blue] {};
\node[draw,inner sep=2.5mm, label=below:,fit=(fwmap2)   (fwmap),color=blue] {};
\end{tikzpicture}

\end{center}

\begin{defn}[Identity optic]
	Given $(A,A')$, the identity optic $\id_{A,A'}^{\mathrm{O}}\colon{}(A,A')\to (A,A')$ is $(1,\lambda^{-1}_A,\lambda_{A'})$:
    we take as internal state $1$, the monoidal unit, then the \revision{first} map is $\lambda^{-1}_A\colon{}A\to 1\otimes A$ and the \revision{second} map is $\lambda_{A'}\colon{}1\otimes A' \to A'$, the left unitors.
\end{defn}

\begin{defn}[Category of optics]
The category $\mathrm{Optic}({\scr{C}})$ has as objects pairs of objects from $\scr{C}$, $(A, A')$, $(B, B'),\ldots$, and as morphisms optics $(M, \fw, \bw)$ from $(A,A')$ to $(B, B')$ as defined in Definition \ref{def:optic}. Sequential composition, the monoidal product and identities are as defined above, and this gives rise to a symmetric monoidal category.

In the case of $\mathrm{Optic}(\Set)$, the monoidal unit is simply $\{*\}\otimes \{*\}$.

\label{def:opticcat}
\end{defn}
See \cite[Theorem 2.0.12]{Riley2018} for a proof that $\mathrm{Optic}({\scr{C}})$ is symmetric monoidal.

\begin{rem}
For our subsequent compositional results for BFFG, we will be working with $\mathrm{Optic}(\Set)$. In this case, since $\Set$ is Cartesian, we have the isomorphism $\mathrm{Optic}(\Set)\cong \mathrm{Lens}(\Set)$ \cite[Proposition~1]{gavranovic2022space}, so we could alternatively have worked with lenses.
\end{rem}

\section{A categorical approach to Backward Filtering Forward Guiding}\label{sec:cat_bffg}

In Section \ref{sec:bffg_intro}, we introduced the BFFG algorithm on a line graph. It consists of two steps:
\begin{itemize}
    \item first, in the \textit{backward pass}, we iteratively compute $g_{i-1 }= \tilde\kappa g_i$;
    \item secondly, in the \textit{forward pass} we iteratively apply the  $h$-transform from Definition \ref{def:htransform} to $(\kappa_i, g_i)$.
\end{itemize}
In this section, we will formalise both steps by defining formally \revision{first} and \revision{second} maps, using the framework of Optics as defined in Section~\ref{sec:optics}. 
The first pass will correspond to the backward pass of BFFG, and the second pass will correspond to the forward pass of BFFG.

\subsection{First and second maps}
In the upcoming definitions we assume we are given Markov kernels $\kappa\colon S \rightarrowtriangle S'$ and $\tilde\kappa\colon S \rightarrowtriangle S'$, where 
$S = (E, \fB)$ and $S' = (E', \fB')$ are measurable spaces. 
Recall that a measure $\mu$ is \textit{s-finite} if there exists a countable family of finite measures $(\mu_n)$ with $\mu = \sum_n \mu_n$.

For a measurable space $S = (E, \fB)$, let $\bM_+(S)$ denote the set of strictly positive measurable functions on $S$. Let $\calM(S)$ denote the set of s-finite measures on $S$. The class of s-finite measures is convenient to work with, since it is closed under \textit{pushforwards} of measurable maps and Markov kernels \cite[Proposition~7]{Staton2017}.

\begin{defn}\label{def:backwardmap} For a Markov kernel $\tilde \kappa\colon S \rightarrowtriangle S'$, define $\backw{\tilde\kappa} \colon \bM_+(S') \to \bM_+( S')\times \bM_+(S)$ by 
\begin{equation}\label{eq:backw}
\backw{\tilde \kappa}(g) = \left(g, \tilde \kappa g\right).
\end{equation}

\end{defn}

Intuitively, we should imagine that the kernel $\tilde\kappa$ is an approximation to the kernel $\kappa$, and then $g$ represents an approximation to the functions $h$ defined through \eqref{eq: hn2}.

This map returns both the  pullback $\tilde \kappa g$  and an appropriate {\it message} $g$ for the map $\forw{\kappa,\tilde\kappa}$ specified in the following definition. This corresponds to the backwards map in BFFG, \revision{and will be the first map, $\fw$,} in terms of the optic (Definition~\ref{def:optic}).
\begin{defn}\label{def:forwardmap}
For Markov kernels $\kappa,\tilde\kappa\colon S \rightarrowtriangle S'$,  define $\forw{\kappa,\tilde\kappa }\colon  \bM_+( S')\times \calM(S)\to \calM(S')$ by 
\begin{equation}\label{eq:forw}
 \forw{\kappa,\tilde\kappa}(g, \mu)  = \nu, \quad \text{with} \quad \nu(\!\dd y) := \int_S \frac{g(y)}{\tilde\kappa g(x)} \mu(\!\dd x) \kappa(x, \dd y). 
 \end{equation}

We note that $\forw{\kappa,\tilde \kappa}$ is well-defined by considering the following decomposition: given $\mu=\sum_n \mu_n$ for finite measures $(\mu_n)$, consider for each $s\in \mathbb N$, $B_s :=\{(x,y)\colon{}s\le m(x,y)<s+1\}\subset S\otimes S'$ with $m(x,y):= \frac{g(y)}{\tilde\kappa g(x)}$; this is well-defined by strict positivity of $g$. Then, $\nu_{n,s}(\dif y) := \int m(x,y) \mu_n(\dif x) \kappa(x,\dif y) 1_{B_s}(x,y)$ is a finite measure, with $\sum_{n,s}\nu_{n,s}=\nu$.

This corresponds to the forwards pass in BFFG, and in terms of the optic will be \revision{the second map, $\bw$,} Definition~\ref{def:optic}.
\end{defn}

The pairing of the \revision{first} and \revision{second} maps can be conveniently summarised by the following diagram, which is analogous to \eqref{eq:optic}:

\begin{center}
\begin{tikzpicture}[style={scale=1.5}]
	\begin{pgfonlayer}{nodelayer}
		\node [style=morphism] (bwmap) at (0, 0) {$\cB_{\tilde\kappa}$};
		\node [style=none] (in_bw) at (-1, 0) {$\bB_1$};	
		\node [style=none] (out_bw) at (1, 0) {$\bB_0$};	
		
		\node [style=morphism] (fwmap) at (0, -1) {$\cF_{\kappa,\tilde\kappa}$};
		\node [style=none] (in_fw) at (-1, -1) {$\bM_1$};	
		\node [style=none] (out_fw) at (1, -1) {$\bM_0$};	

		\node [style=none] (message) at (0.25, -0.5) {$g$};	
	\end{pgfonlayer}	
	
	\begin{pgfonlayer}{edgelayer}
		\draw [style=edge] (in_bw) to (bwmap);
		\draw [style=edge] (bwmap) to (out_bw);		
		
		\draw [style=edge] (fwmap) to (in_fw);
		\draw [style=edge] (out_fw) to (fwmap);		

		\draw [style=edge] (bwmap) to (fwmap) ;		
	\end{pgfonlayer}

\end{tikzpicture}
\end{center}

Note that in the case $\tilde\kappa = \kappa$, if $\mu$ is a probability measure, then $\cF_{\kappa,\kappa}(g,\mu)$ is a  probability measure as well. 
In general, even if $\mu$ is a probability measure, the output of $\forw{\kappa, \tilde\kappa}$ will typically {\it not} be a probability measure. 
In that case, we can reinterpret $\forw{\kappa,\tilde\kappa}$ as producing a {\it weighted} probability measure: 
if $\varpi\ge 0$ and $\mu$ is a probability measure, then
\[
 \forw{\kappa,\tilde\kappa}(g, \varpi \cdot \mu)(\!\dd y) = (\varpi w_{\kappa,\tilde\kappa}(g, \mu)) \cdot \nu(\!\dd y), \]
where the {\it weight} $w_{\kappa,\tilde\kappa}(g, \mu)$ and probability measure $\nu$ are defined by 
\begin{equation}\label{eq:wnu}
\begin{split}
 w_{\kappa,\tilde\kappa}(g,\mu) &= \int_S \int_{S'}  \frac{g(y)}{\tilde\kappa g(x)} \kappa(x,\dd y) \mu(\dif x) =  \int_S \frac{(\kappa g)(x)}{(\tilde \kappa g)(x)} \mu(\dif x)\\
\nu(\!\dd y) &=  w^{-1}_{\kappa,\tilde\kappa}(g, \mu) \int_S  \frac{g(y) }{(\tilde \kappa g)(x)  }  \mu(\!\dd x) \kappa(x, \dd y), 
 \end{split}
 \end{equation}
 provided we have $0<w_\kappa(m,\mu)<\infty$.
Hence, by using an approximate $\tilde\kappa$ in $\backw{\tilde\kappa}$, we ``pick up'' a weight $w_{\kappa,\tilde\kappa}(g,\mu)$.

\begin{defn}
    Given kernels $\tilde\kappa, \kappa\colon{}S\rightarrowtriangle S'$, our triple $(\bM_+(S'),\backw{\tilde\kappa}, \forw{\kappa,\tilde\kappa})$ defines an optic $(\bM_+(S'),\calM(S'))\to (\bM_+(S),\calM(S))$ in the category $\mathrm{Optic}(\Set)$.
This optic will be denoted by $\scr{O}(\kappa, \tilde\kappa)$.  
\label{def:scrO}
\end{defn} 

\begin{rem}
The BFFG algorithm (as in Algorithm~\ref{alg:bffg}) is related to our above construction as follows:
    \begin{itemize}
        \item Given $g_{n+1}=\tilde p_{n+1}(\cdot, x_{n+1})$ analogous to \eqref{eq: hn2}, the backward filtering step of Algorithm~\ref{alg:bffg}, Line~17 is equivalent to the composition of first maps $\backw{\tilde\kappa_{1}}, \dots,\backw{\tilde\kappa_{n}}(g_{n+1})$, using the composition rule \eqref{eq:fw_comp}. This produces a sequence of guiding functions $(g_0,\dots, g_n)$ needed for the second pass. 
        \item The forward guiding step of Algorithm~\ref{alg:bffg}, Line~18 is equivalent to sampling from the compositions of second maps $\cF_{\kappa_n,\tilde\kappa_n},\dots,\cF_{\kappa_1,\tilde\kappa_1}(\delta_{\hat x_0})$ using the composition rule \eqref{eq:bw_comp}, given the output $(g_0,\dots,g_n)$ of the first maps.
    \end{itemize}
    Thus the overall BFFG algorithm is precisely given by
    composing the optics $(\bM_+(S_i),\backw{\tilde\kappa_i}, \forw{\kappa,\tilde\kappa_i})_{i=1}^n$ as given in Definition~\ref{def:scrO}.  
\end{rem}

\begin{defn}
    Since we will be working with pairs of kernels $\kappa, \widetilde\kappa$, we define the category $\mathrm{Stoch}_2$ to be the category whose objects are measurable spaces, as in $\Stoch$, but whose morphisms are now \textit{ordered pairs} of kernels $(\kappa, \widetilde\kappa)\colon{}S\rightarrowtriangle S'$ with the same source and target. Composition is then performed component-wise: $(\kappa, \widetilde\kappa)\circ (\kappa', \widetilde\kappa')=(\kappa\circ \kappa', \widetilde\kappa\circ \widetilde \kappa ')$, {and similarly the monoidal product is also defined pairwise on each component}. The required associativity and identities are directly inherited from $\Stoch$.
\end{defn}


\subsection{Functoriality of $\scr{O}$: equivalence of sequential composition}\label{sec:seq}

It turns out that the BFFG algorithm preserves sequential composition from kernels to optics. Formally, we have the following:

\begin{thm}\label{thm:functor}
    The BFFG algorithm defines a contravariant functor $\scr{O}\colon{} (\Stoch_2)^{\mathrm{Op}}\to \Optic(\Set)$:
    \begin{itemize}
        \item on objects, given a measurable space $S$, we have $\scr{O}(S)= \revision{(\bM_+(S),\calM(S))}$;
        \item on arrows, given $(\kappa, \tilde\kappa)\colon{}S\rightarrowtriangle S'$, we have $\scr{O}(\kappa,\tilde\kappa)\colon{}\scr{O}(S')\to \scr{O}(S)$ is as defined in Definition~\ref{def:scrO}.
    \end{itemize}
    \label{prop:seq}
\end{thm}

To spell this out, consider the composition of two Markov kernels, say $\kappa_{0,1}$ and $\kappa_{1,2}$, and corresponding kernels $\tilde\kappa_{0,1}$ and $\tilde\kappa_{1,2}$. We introduce the short-hand notation $\cB_{i,i+1}=\cB_{\tilde\kappa_{i,i+1}}$ and  $\cF_{i,i+1}=\cF_{\kappa_{i,i+1},\tilde \kappa_{i,i+1}}$ for $i\in \{0,1\}$.
The optics $\scr{O}(\kappa_{0,1}, \tilde\kappa_{0,1})$ and $\scr{O}(\kappa_{1,2}, \tilde\kappa_{1,2})$ induced by $(\tilde\kappa_{0,1}, \kappa_{0,1})$ and $(\tilde\kappa_{1,2}, \kappa_{1,2})$ respectively can be composed as follows:

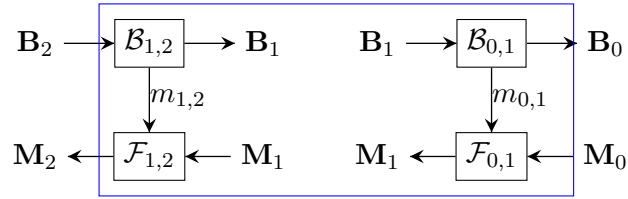
\begin{center}

\begin{tikzpicture}[style={scale=1.5}]
	\begin{pgfonlayer}{nodelayer}
		\node [style=morphism] (bwmap2) at (0.5, 0) {$\cB_{1,2}$};
		\node [style=none] (in_bw2) at (-1, 0) {$\bB_2$};	
		
		\node [style=morphism] (fwmap2) at (0.5, -1.5) {$\cF_{1,2}$};
		\node [style=none] (in_fw2) at (-1, -1.5) {$\bM_2$};	
		\node [style=none] (message2) at (0.75, -0.75) {$g_{1,2}$};	
	\end{pgfonlayer}	
	\begin{pgfonlayer}{edgelayer}
		\draw [style=edge] (in_bw2) to (bwmap2);
		
		\draw [style=edge] (fwmap2) to (in_fw2);
		\draw [style=edge] (bwmap2) to (fwmap2);		
	\end{pgfonlayer}
	\begin{pgfonlayer}{nodelayer}
		\node [style=morphism] (bwmap) at (2.5, 0) {$\cB_{0,1}$};
		\node [style=none] (out_bw) at (4, 0) {$\bB_0$};	
		
		\node [style=morphism] (fwmap) at (2.5, -1.5) {$\cF_{0,1}$};
		\node [style=none] (out_fw) at (4, -1.5) {$\bM_0$};	
		\node [style=none] (message) at (2.75, -0.75) {$g_{0,1}$};	
	\end{pgfonlayer}	
	
	\begin{pgfonlayer}{edgelayer}
		\draw [style=edge] (bwmap2) to (bwmap);
		
		\draw [style=edge] (fwmap) to (fwmap2);
		
		\draw [style=edge] (bwmap) to (out_bw);		
		\draw [style=edge] (out_fw) to (fwmap);		
		\draw [style=edge] (bwmap) to (fwmap);		
	\end{pgfonlayer}

	\node [style=none] (B_label) at (1.5, 0.2) {$\bB_1$};
	\node [style=none] (Bprime_label) at (1.5, -1.3) {$\bM_1$};

    \node[draw,inner sep=2.5mm, label=below:,fit=(bwmap)  (bwmap2),color=blue] {};
    \node[draw,inner sep=2.5mm, label=below:,fit= (fwmap) (fwmap2) ,color=blue] {};
    
\end{tikzpicture}

\end{center}

However, since Markov kernels can be composed -- using the Chapman-Kolmogorov equation \eqref{eq:chapman} -- there is an alternative way to compose. Defining the kernels $\cB_{0,2}=\cB_{\tilde \kappa_{0,2}}$, $\cF_{0,2}=\cF_{\kappa_{0,2},\tilde \kappa_{0,2}}$ with composed kernels $\kappa_{0,2}=\kappa_{0,1}\kappa_{1,2}$ and  $\tilde \kappa_{0,2}=\tilde \kappa_{0,1}\tilde \kappa_{1,2}$ , we have
\begin{center}
\begin{tikzpicture}[style={scale=1.5}]
	\begin{pgfonlayer}{nodelayer}
		\node [style=morphism] (bwmap) at (0, 0) {$\cB_{0,2}$};
		\node [style=none] (in_bw) at (-1, 0) {$\bB_2$};	
		\node [style=none] (out_bw) at (1, 0) {$\bB_0$};	
		
		\node [style=morphism] (fwmap) at (0, -1) {$\cF_{0,2}$};
		\node [style=none] (in_fw) at (-1, -1) {$\bM_2$};	
		\node [style=none] (out_fw) at (1, -1) {$\bM_0$};	

		\node [style=none] (message) at (0.25, -0.5) {$g_{0,2}$};	
	\end{pgfonlayer}	
	
	\begin{pgfonlayer}{edgelayer}
		\draw [style=edge] (in_bw) to (bwmap);
		\draw [style=edge] (bwmap) to (out_bw);		
		
		\draw [style=edge] (fwmap) to (in_fw);
		\draw [style=edge] (out_fw) to (fwmap);		

		\draw [style=edge] (bwmap) to (fwmap) ;		
	\end{pgfonlayer}

\end{tikzpicture}.
\end{center}

Hence, there are two ways to compose:
\begin{enumerate}
\item {\bf construction 1:} compose Markov kernels to get $\kappa_{0,2}=\kappa_{0,1}\kappa_{1,2}$ (similarly $\tilde\kappa_{0,2}=\tilde\kappa_{0,1} \tilde\kappa_{1,2}$), then form the optic;
	\item {\bf construction 2:} compose optics $(\bM_+(S'), \cB_{0,1}, \cF_{0,1})$ and $(\bM_+(S''), \cB_{1,2}, \cF_{1,2})$.
\end{enumerate}

The content of Theorem~\ref{prop:seq} is that these two constructions give rise to the same optic. The proof is broken into two steps.

\begin{lem}
    Both optics are in the same equivalence class. That is
	\[ \scr{O}(\kappa_{0,2}, \tilde\kappa_{0,2}) \sim  \scr{O}(\kappa_{1,2}, \tilde\kappa_{1,2})\circ \scr{O}(\kappa_{0,1}, \tilde\kappa_{0,1}), \]
the right-hand-side denoting composition of optics. 	
\end{lem}

\begin{proof}
For the construction 1, we see that it has internal state $M:= \bM_+(S'')$, with
\begin{equation*}
    \backw{0,2}(g) = \left( g,\tilde \kappa_{0,2}g \right),\quad g\in \bM_+(S''),
\end{equation*}
and for $g\in \bM_+(S'')$ and $\mu\in \calM(S)$,
\begin{equation*}
\begin{split}
    \forw{0,2} (g,\mu) &=\nu,\\
    \nu(\dif y) &= \int_S \frac{g(y)}{\tilde\kappa_{0,2}g(x)} \mu(\dif x) \kappa_{0,2}(x,\dif y).
\end{split}
\end{equation*}

On the other hand, with construction 2 by composing optics we find we have inner state $N:= \bM_+(S'')\otimes \bM_+(S') $ and maps
\begin{equation*}
    \backw{0,2}'\colon{} g\mapsto \left( g,\tilde\kappa_{1,2}g,\tilde\kappa_{0,2}g \right),\quad g\in \bM_+(S''),
\end{equation*}
and for $(g,g',\mu)\in \bM_+(S'')\otimes \bM_+(S')\otimes \calM(S)$,
\begin{equation}
\begin{split}
    \forw{0,2}'\colon{}(g,g',\mu) &\mapsto \nu',\\
    \nu'(\dif y) &= \int_S \frac{\mu(\dif x)}{\tilde\kappa_{0,1}g'(x)}\int_{S'}\kappa_{0,1}(x,\dif y')\cdot  \frac{g'(y')}{\tilde \kappa_{1,2}g(y')} \cdot \kappa_{1,2}(y',\dif y) g(y).
\end{split}
\label{eq:nuprime}
\end{equation}

To show that these are in the same equivalence class, we define the residual morphism $r\colon{}M\to N$ via
\begin{equation*}
    r\colon{} \bM_+(S'') \ni g \mapsto \left(g, \tilde \kappa_{1,2}g \right)\in \bM_+(S'')\otimes \bM_+(S').
\end{equation*}
With this choice of $r$, it is easy to see that the required commutation relations on $r$ in Definition~\ref{def:optic} do hold; for instance in the definition of $\nu'$ in \eqref{eq:nuprime}, given $g'=\tilde\kappa_{1,2}g$, the inner fraction simplifies to $1$ and we see $\nu'=\nu$ in this case.
\end{proof}

In order to complete the proof of Theorem~\ref{prop:seq}, we also need to establish preservation of identities.
\begin{lem}
    For any measure space $S$, writing $\id_S$ for the identity kernel $x\mapsto \delta_x$ in $\Stoch$, we have that $\scr{O}(\id_S, \id_S)$ is the identity optic $\id^{\mathrm O}_{\scr{O}(S)}$ in $\Optic(\Set)$.
\end{lem}
\begin{proof}
    For the \revision{first} pass, we have $\backw{\id_S}(g)=(g,g)$.
    For the \revision{second} pass, we have $\forw{\id_S,\id_S}(g,\mu)=\nu$ with $\nu(\dif x)=\mu(\dif x)$.

    By trivially discarding the message $g$, that is, $r\colon{}\bM_+(S)\to 1=\{*\}$ maps $g\mapsto *$, we see this corresponds to the identity optic on $(\bM_+(S),\calM(S))$.
\end{proof}

\subsection{Equivalence of parallel composition}\label{sec:par}
In the previous section we established that BFFG is functorial, namely it preserves sequential composition. 
We now show that BFFS (when $\tilde\kappa=\kappa$) is a \textit{lax monoidal functor}: it preserves parallel composition.
Intuitively, lax monoidality of BFFS is related to the fact that conditional independence is preserved under conditioning on the leaves. We would not expect more than lax monoidality of BFFS (for instance, bilaxity)  due to the phenomenon of confounding: conditioning on joint observations introduces dependence between formerly independent states.
The failure of lax monoidality of BFFG is related to the additional need to track global weights; see Remark~\ref{rem:lax_fail}.

Recall \eqref{eq:markovpar} defining a tensor product of Markov kernels. 
Intuitively, there are two constructions involving the `product' of two Markov kernels $\kappa_1$ and $\kappa_2$
\begin{enumerate}
	\item {\bf construction 1:} take the product of Markov kernels to get $\kappa_{1\otimes 2}=\kappa_{1} \otimes \kappa_{2}$ (similarly for $\tilde\kappa_{1\otimes 2}=\tilde\kappa_{1} \otimes \tilde\kappa_{2}$), then form the corresponding optic;
	\item {\bf construction 2:} take the tensor product of the optics $(M_1, \cB_{1} ,\cF_{1})$ and $(M_{2}, \cB_{2}, \cF_{2})$.
\end{enumerate}

The goal of this section is to prove our following main result. The crucial restriction here is that we have only a single kernel $\kappa=\tilde \kappa$, without an additional guiding kernel: in this setting, the category $\Stoch_2$ can be identified with $\Stoch$.

\begin{thm}\label{thm:lax}
    In the setting where $\tilde\kappa = \kappa$, the BFFG algorithm defines a lax monoidal functor  $\scr{O}\colon{} (\Stoch)^{\mathrm{Op}}\to \Optic(\Set)$.
    \label{thm:lax_mf}
\end{thm}
\begin{proof}
    Since we have already established in Theorem~\ref{prop:seq} that $\scr{O}\colon{} (\Stoch_2)^{\mathrm{Op}}\to \Optic(\Set)$ defines a functor, what remains to establish lax monoidal functoriality is defining an appropriate morphism $\epsilon$ and natural transformation $\mu$, which together satisfy associativity and unitality.

    Recall that the monoidal unit in $\Stoch$ is $1_\Stoch = (\{*\}, (\emptyset , \{*\}))$ and in $\Optic(\Set)$ we have $1_{\Optic(\Set)} = \{*\}\otimes\{*\}$.

For each pair of objects $S,T$ in $\Stoch$, define the optic 
$${\mu_{S,T}}\colon{}\scrO(S)\otimes\scrO(T)\to\scrO\left(S\otimes T\right)$$ as follows: we take trivial internal state, and the following \revision{first} and \revision{second} maps:
\begin{equation*}
    \begin{split}
        \fw\colon{} \bM_+(S)\otimes \bM_+(T)\to \bM_+(S\otimes T),\quad \left(g(s),h(t)\right)&\mapsto g(s)\cdot h(t)\\
        \bw\colon{} \calM(S\otimes T)\to \calM(S)\otimes \calM(T),\quad \bar \mu(\dif s,\dif t) &\mapsto \left(\bar \mu (\dif s, T), \bar\mu (S,\dif t)\right),
    \end{split}
\end{equation*}
noting that the marginal measures of an $s$-finite measure remain $s$-finite.

In other words, a pair of functions on the individual spaces is mapped to the pointwise product function on the joint space, and a measure on the joint space is mapped to its two marginals on the individual spaces.

In order to prove naturality of $\mu$, we need to establish the following naturality square: for kernels $\kappa_1, \tilde\kappa_1 \colon{} S\rightarrowtriangle S'$ and $\kappa_2, \tilde\kappa_2 \colon{} T\rightarrowtriangle T'$ ,

\begin{equation}\label{eq:naturality_sq}
\begin{tikzcd}[column sep = 16ex]
    \scrO(S')\otimes\scrO(T') \arrow[r, "\scrO(\kappa_1;\tilde\kappa_1) \otimes {\scrO(\kappa_2;\tilde\kappa_2) } "] \arrow[d, "\mu_{S',T'}"] 
    & \scrO(S)\otimes\scrO(T) \arrow[d, "\mu_{S,T}"] \\
    \scrO(S'\otimes T')\arrow[r, "\scrO \left(\kappa_1\otimes \kappa_2;\tilde\kappa_1\otimes \tilde\kappa_2\right)"]
    & \scrO(S\otimes T)  
\end{tikzcd}.
\end{equation}
In Lemma~\ref{lem:naturality} below, we demonstrate that this naturality square does hold in the assumed case when $\tilde\kappa_1=\kappa_1$ and $\tilde\kappa_2=\kappa_2$.

Thus in order to complete the proof 
 of Theorem~\ref{thm:lax_mf}, we need to  verify the associativity and unitality conditions of Definition~\ref{defn:laxmf}.

To this end, define the morphism $\epsilon\colon{} 1_{\Optic(\Set)}\to \scr{O}(1_\Stoch)=\bM_+(\{*\})\otimes \calM(\{*\})$ to be the optic with trivial internal state and the following \revision{first} and \revision{second} maps:
\begin{equation*}
    \begin{split}
        \fw\colon{} \{*\}\to \bM_+(\{*\}),\quad *\mapsto 1\\
        \bw\colon{} \calM(\{*\})\to \{*\},\quad \mu \mapsto *,
    \end{split}
\end{equation*}
that is, \revision{in the first pass} $*$ is mapped to the constant function 1 and \revision{in the second pass} the measure is simply deleted.

In Lemma~\ref{lem:associativity}, we prove associativity of $\mu$. In Lemma~\ref{lem:unitality}, we prove that the unitality conditions on $\epsilon,\mu$ hold. 

Thus the proof of Theorem~\ref{thm:lax} is complete.
\end{proof}

\begin{lem}\label{lem:naturality}
    In the case when $\tilde\kappa_1=\kappa_1$ and $\tilde\kappa_2=\kappa_2$, our family of morphisms $\mu_{S,T}$ defines a natural transformation $\scrO(\cdot)\otimes \scrO(\cdot) \to \scrO(\cdot \otimes \cdot)$.
\end{lem}
\begin{proof}
    
Let's compute the two sides of the naturality square \eqref{eq:naturality_sq}:
\begin{itemize}
    \item $(\downarrow, \rightarrow)$: here we first map using the component $\mu_{S',T'}$, and then use the BFFG functor on a pair of product kernels. This composite optic has the following components:
    \begin{equation}
        \begin{split}
            \text{internal state: } N &=\bM_+(S'\otimes T') ;\\
            \fw: \bM_+(S')\otimes \bM_+(T')&\to N\otimes \bM_+(S\otimes T),\\
            (g_1,g_2) &\mapsto \left( g_1 \cdot g_2,\tilde \kappa_1 g_1 \cdot \tilde\kappa_2 g_2 \right);\\
            \bw: N\otimes \calM(S\otimes T)&\to \calM(S')\otimes \calM(T'),\\   
            \left( \bar g(s',t'),\bar \mu(\dif s,\dif t) \right) &\mapsto 
            \left (\nu(\dif s', T'), \nu(S', \dif t') \right), \\
            \text{where } \nu(\dif s', \dif t') =  \int_{S\otimes T} &\frac{\bar{g}(s',t')}{\left[\left(\tilde{\kappa}_{1}\otimes\tilde{\kappa}_{2}\right)\bar{g}\right](s,t)}\bar{\mu}(\dif s,\dif t)\kappa_{1}(s,\dif s')\kappa_{2}(t,\dif t').
        \end{split}
        \label{eq:lax1}
        \end{equation}
{The distributions \eqref{eq:lax1} can be naturally interpreted as the marginal distributions of the output of the guiding pass (the second map).}  
    
    \item $(\rightarrow,\downarrow)$: here we first form the product of two individual optics formed from BFFG, and then use the component $\mu_{S,T}$. We have the following components:
    \begin{equation}
    \begin{split}
        \text{internal state: } M &=\bM_+(S')\otimes\bM_+( T') ;\\
        \fw:\bM_+(S')\otimes \bM_+(T')&\to M\otimes \bM_+(S\otimes T),\\
            (g_1,g_2) &\mapsto \left( g_1 , g_2,\tilde \kappa_1 g_1 \cdot \tilde\kappa_2 g_2 \right);\\            
            \bw: M\otimes \calM(S\otimes T)&\to \calM(S')\otimes \calM(T'),\\
            \left(  g_1,g_2,\bar \mu(\dif s,\dif t) \right) \mapsto \bigg( \int_{S}\frac{{g_1}(s')}{\tilde{\kappa}_{1}{g_1}(s)}\bar{\mu}(\dif s,T)&\kappa_{1}(s,\dif s') ,\int_{T}\frac{{g_2}(t')}{\tilde{\kappa}_{2}{g_2}(t)}\bar{\mu}(S,\dif t)\kappa_{2}(t,\dif t')\bigg).
    \end{split}\label{eq:rd_optic}
    \end{equation}
\end{itemize}
In order to argue that these two optics are equal, we need to present a residual morphism $r\colon{}M\to N$ between the internal states, such that the two triangles from Definition~\ref{def:optic} commute. We choose
\begin{equation*}
\begin{split}
    r\colon{}\bM_+(S')\otimes\bM_+( T')&\to \bM_+(S'\otimes T'),\\
    \left(g_{1}(s'),g_{2}(t')\right)&\mapsto g_{1}(s')\cdot g_{2}(t').
\end{split}
\end{equation*}
With this definition, the left-hand triangle of Definition~\ref{def:optic} is immediately seen to hold.

For the right-hand triangle, it is equivalent to requiring that the following is equal to the final expression of \eqref{eq:rd_optic} for any $g_1,g_2, \bar \mu$:
\begin{equation}
         \bigg( \int_{S\otimes T}\frac{{g_1}(s')}{\tilde{\kappa}_{1}{g_1}(s)}\bar{\mu}(\dif s,\dif t)\kappa_{1}(s,\dif s') \,\mathrm{w}_2(t),
            \int_{S\otimes T}\frac{{g_2}(t')}{\tilde{\kappa}_{2}{g_2}(t)}\bar{\mu}(\dif s,\dif t)\kappa_{2}(t,\dif t')\, \mathrm w_1(s)\bigg),
            \label{eq:lax2}
\end{equation}
where
\begin{equation*}
    \mathrm w_i(t) := \frac{\kappa_i g_i (t)}{\tilde \kappa_i g_i(t)}, \quad i=1,2.
\end{equation*}
Naturality thus holds when $\tilde\kappa_1=\kappa_1$ and $\tilde\kappa_2 = \kappa_2$. 
\end{proof}

\begin{rem}
    {We can see in the proof of Lemma~\ref{lem:naturality} why BFFG fails to be a lax monoidal functor in general: without the restriction $\tilde\kappa=\kappa$, the `natural' choice of $\mu$ in fact {fails to be a natural transformation}: the $\mathrm w_i$ appearing in \eqref{eq:lax2} can be interpreted as the importance weights picked up from the \textit{other} branch of the computation, which in general, do \textit{not equal} $1$.
While we do not formally prove that there do not exist any natural transformations that can make BFFG lax monoidal, we conjecture that no such transformations exist.}
\end{rem}

\begin{lem}\label{lem:associativity}
    Our morphisms $\mu$ satisfy the associativity condition.
\end{lem}
\begin{proof}
    We need to check that the associativity hexagon in Definition~\ref{defn:laxmf} holds.

    This is somewhat tedious to verify, but it does hold: we need to check the actions on functions and measures, checking that for these mappings we can take either direction around the hexagon. 

    For functions, since pointwise products of functions is associative, both directions will map a triplet $(f(r),g(s),h(t))$ of functions to the pointwise product $(r,s,t)\mapsto f(r)\cdot g(s)\cdot h(t)$.

    For measures, the marginalisation operation is also associative, in the sense that given a measure on a triple joint space $\bar \mu (\dif r, \dif s, \dif t)$, it does not matter in which order one performs the marginalisation: one can first form the intermediate pair of marginals $\left(  \bar \mu (\dif r, S,T), \bar \mu (R,\dif s, \dif t)\right)$ or the alternative pair $\left(  \bar \mu (\dif r, \dif s,T), \bar \mu (R,S, \dif t)\right)$.
    Given either of these intermediate pairs, further marginalisation will result in the same set of marginals, namely $\left(  \bar \mu (\dif r, S,T), \bar \mu (R,\dif s, T),\mu (R,S, \dif T)\right)$.
\end{proof}

\begin{lem}\label{lem:unitality}
    Our morphisms $\epsilon, \mu$ satisfy the unitality conditions.
\end{lem}
\begin{proof}
    We check now that the two unitality squares of Definition~\ref{defn:laxmf} hold.

    For the left square, for functions we are comparing the identity mapping on a function $g$ to the product $g\cdot 1$ with the constant 1 function. On measures, given a measure $\mu(\dif s)$, we first marginalise to obtain $\left( \mu(S)\cdot \delta_*, \mu \right)$ and then delete, resulting in $(*,\mu)$, as required.

    The right square is essentially identical.
\end{proof}

\begin{rem}
    The previous two lemmas are valid in the general case of $\Stoch_2$ and do not require the restriction to $\tilde\kappa=\kappa$.
\end{rem}

\begin{rem}\label{rem:failure lax}
Failure of lax-monoidality is not surprising. Clearly, upon first extracting marginals from the joint measure $\bar\mu$ and then applying BFFG on each of the marginals, one loses the dependency structure of $\bar \mu$, which is retained when applying BFFG on the joint space. However, when we have a product measure, $\bar\mu(\dif s, \dif t)=\bar\mu(\dif s, T) \bar\mu(S, \dif t)$, Equation \eqref{eq:lax2} reads
\[
      \bigg( \int_{S}\frac{{g_1}(s')}{\tilde{\kappa}_{1}{g_1}(s)}\bar{\mu}(\dif s,T)\kappa_{1}(s,\dif s') C_2 ,\int_{T}\frac{{g_2}(t')}{\tilde{\kappa}_{2}{g_2}(t)}\bar{\mu}(S,\dif t)\kappa_{2}(t,\dif t')C_1\bigg)
\]
with $C_1$ and $C_2$ constants given by
\[ C_1=\int_S \mathrm w_1(s) \bar\mu(\dif s, T) \quad \text{and}\quad C_2=\int_T \mathrm w_2(t) \bar\mu(S, \dif t). \]
This is to be compared to the expression for $\bw$ in \eqref{eq:rd_optic}. 
Algorithmically, since in applications of BFFG we typically initialise from a Dirac mass -- which is indeed a product measure $\delta_{(x,y)}=\delta_x\otimes \delta _y$ -- these constants are the only differences between the two sides of the naturality square. The output of the corresponding implementations of the BFFG algorithm is in both cases a sample from $x$ from $P^\circ$ its likelihood, up to a constant factor. 
Considering the form of the final estimator $T$  \eqref{eq:selfnormalised_is}, these constants will ultimately cancel out, and hence are irrelevant. Thus, while lax-monoidality fails from a mathematical perspective, from a statistical perspective this failure is not an issue, and we can consider either implementation of BFFG. 
    \label{rem:lax_fail}
\end{rem}

\begin{rem}
We comment on links with the recent works of \cite{Smithe2020, Smithe2021, Braithwaite2023}, which {contain} parallels of our results, when they are specialised to the BFFS setting. In these cited works, the general procedure of Bayesian inference (namely obtaining a posterior from a prior and data) is placed on a general categorical footing, and shown to have an optical structure, similar to our treatment of BFFG. Indeed \cite[Theorem 5.2]{Smithe2020} also demonstrates functoriality of Bayesian inversion into a category of optics; our Theorem~\ref{thm:functor} in the case of BFFS can thus be seen as a parallel of this result. For example, in the case of BFFS ($\tilde\kappa=\kappa$), we can dually interpret the function $g$ as an (improper) `prior' density, dually interpret the measure $\mu$ as a `test function', and then the measure $\forw{\kappa}(g,\mu)$ from the second pass can be interpreted as a Bayesian inversion of the channel $\kappa$ along `state' $g$. Lax {monoidality} of the Bayesian inversion functor has also been established in \cite[Remark~22]{Braithwaite2023}, under the assumption of the existence of \textit{support objects} (for instance, in the case of finite state spaces or Gaussian kernels). This result, when applicable, thus parallels our Theorem~\ref{thm:lax_mf}, which we have proven in the setting of $s$-finite measures and measurable functions.
\end{rem}

\section{Examples}\label{subsec:examples}

The results in this paper show that we may split an edge into multiple edges or, conversely, collapse multiple edges: the BFFG-algorithm is not affected by this. Below we briefly describe some examples that illustrate the relevance of this property.

\begin{figure}
\begin{center}
\begin{tikzpicture}[style={scale=0.52}]
	\tikzstyle{empty}=[fill=white, draw=black, shape=circle,inner sep=1pt, line width=0.7pt]
	\tikzstyle{solid}=[fill=black, draw=black, shape=circle,inner sep=1pt,line width=0.7pt]
	\begin{pgfonlayer}{nodelayer}
		\node [] (00-) at (-8, 0) {};
		\node [style=solid,label=below:{$t_{i-1}$},] (00) at (-6, 0) {};
		\node [style=empty,label={$v_{i-1}$},] (0obs) at (-6, 1.5) {};

		\node [style=solid,label=below:{$t_i$},] (0) at (-4, 0) {};

		\node [style=solid,label=below:{$t_{i+1}$}] (1) at (-2, 0) {};
		\node [style=empty,label={$v_{i+1}$},] (2obs) at (-2, 1.5) {};

        \node [] (00-a) at (2,0){};
        \node [style=solid,label=below:{$t_{i-1}$},] (00a) at (4, 0) {};
        \node [style=solid, label=below:{$t_{i+1}$}] (1a) at (7,0) {};
        \node [] (2a) at (9,0) {};

        \node [style=empty,label={$v_{i-1}$},] (0aobs) at (4, 1.5) {};
        \node [style=empty,label={$v_{i-1}$},] (1aobs) at (7, 1.5) {};

        \node [label=above:$\kappa_{i-1}$] at (-5,0){};
        \node [label=above:$\kappa_{i}$] at (-3,0){};

        \node [label=above:$\kappa_{i-1}\kappa_i$] at (5.5,0){};
    
	\end{pgfonlayer}
	\begin{pgfonlayer}{edgelayer}
		\draw [style=edge] (00-) to (00);
		\draw [style=edge] (00) to (0) ;
		\draw [style=edge] (0) to (1);
		\draw [style=edge] (1) to (2);

		\draw [style=edge,color=blue] (00) to (0obs);
		\draw [style=edge,color=blue] (1) to (2obs);

        \draw [style=edge] (00-a) to (00a);
        \draw [style=edge] (00a) to (1a);
        \draw [style=edge] (1a) to (2a);

        \draw [style=edge,color=blue] (00a) to (0aobs);
		\draw [style=edge,color=blue] (1a) to (1aobs);
	\end{pgfonlayer}
\end{tikzpicture}

\caption{Visualisation for Examples~\ref{ex:missing} and \ref{ex:intract}. On the left, there are two separate kernels $\kappa_{i-1}:S_{i-1}\to S_i$ and $\kappa_i:S_i\to S_{i+1}$, with no observation at point $t_i$. On the right, there is a single (composed) kernel $\kappa_{i-1}\kappa_i:S_{i-1}\to S_{i+1}$.}
\label{fig:missing}
\end{center}
\end{figure}

\begin{ex}[Missing observation]
Consider the hidden Markov model setting and suppose that at a hidden node $t_i$ an observation is not recorded,  perhaps due to a failing sensor -- this situation is depicted on the left in Figure~\ref{fig:missing}, where the observation $v_i$ is absent. In this case, Theorem~\ref{thm:functor} tells us that can we can collapse the edges $(t_{i-1}, t_i)$ and $(t_i, t_{i+1})$ to $(t_{i-1}, t_{i+1})$ 
in the latent process, corresponding to the right diagram of Figure~\ref{fig:missing}. Specifically, instead of performing two separate steps of BFFG with kernels $\kappa_{i-1}$ and $\kappa_i$ (and composing optics), we could instead perform a single step of BFFG with the composite kernel $\kappa_{i-1}\kappa_i$.
Theorem \ref{thm:functor} ensures that the overall algorithm is not affected by this. 
\label{ex:missing}
\end{ex}

\begin{ex}[Intractable composition]

   Consider the setting where a transition over an edge is captured by the kernel $\kappa$, which is intractable, but can be expressed as  $\kappa=\kappa_1\kappa_2$ where both $\kappa_1$ and $\kappa_2$ are tractable.  This suggest splitting the edge into two edges, where the added vertex in the middle represents  a latent variable. Writing $\phi(\mu, \sigma)$ for a Gaussian random variable with mean $\mu$ and covariance $\sigma$, a concrete example is $\kappa_1(x,\dd y)=\phi(\dif y; \alpha x, \sigma^2)$ and $\kappa_2(x,\dd y)=\phi(\dif y;\mu(x), \nu^2)\dd y$.
   In this case
   \[ {\kappa(x,\dd y)} =\int \kappa_1(x, \dd z) \kappa_2(z,\dd y) = \int \phi(\dif z; \alpha x, \sigma^2) \phi(\dif y; \mu(x), \nu^2) .\]
   Only for very specific choices of the map $\mu$ this can be computed in closed form, hampering the backward step. However, due to the structure of $\kappa$, we can split the backward step into the composition of two backward steps, where we take  $\tilde\kappa_1=\kappa_1$ and $\tilde\kappa_2(x,\dd y) =\phi(\dif y;\tilde\mu, \sigma^2)$ for a user-specified value of $\tilde\mu$. Only the second of these steps involves an approximation.
Diagrammatically, this corresponds with beginning from the \textit{right} diagram of Figure~\ref{fig:missing}, and introducing an auxiliary intermediary state, corresponding to the left diagram of Figure~\ref{fig:missing}.
   Again, Theorem \ref{thm:functor} allows us to conclude that BFFG remains the same when introducing this latent variable, performing two steps of BFFG with kernels $\kappa_1,\kappa_2$ respectively and then composing the optics.
   \label{ex:intract}
\end{ex}

\begin{ex}[SDE example]
Consider an edge $e$. On this edge, the transition is governed by the evolution of the stochastic differential equation (SDE) $\dd X_t =\mu(X_t) \dd t + \dd W_t$ for $t \in [0,\tau_e]$.
Note that in general the transition densities for this SDE are not known in closed form. However, on short time intervals, an Euler discretisation provides an accurate approximation. As a result, one can include additional vertices that represent times $i\tau_e/n$, $i=1,\ldots, n-1$, for some (large) specified value of $n$.  Similar to the previous example, the backward filtering step can now be split into $n$ steps, where the approximate kernels $\tilde\kappa$ on vertices correspond to the Euler-discretised SDE. 
\end{ex}

\appendix
\section{Categorical definitions}
\label{app:cat}

We start with the definition of a category, as in \cite[(Definition 1.1.1)]{pierce1991basic}.
\begin{defn}
	A {\bf category} $\mathcal C$ comprises of
	\begin{enumerate}
		\item a collection of objects;
		\item a collection of arrows (also called morphisms);
		\item operations assigning to each arrow $f$ an object $\mathrm{dom} f$, its domain, and an object $\mathrm{cod} f$, its codomain, succinctly written  $f\colon A\to B$ where $A=\mathrm{dom} f,B=\mathrm{cod} f$;
		\item a composition operator assigning to each pair of arrows $f$ and $g$ with $\mathrm{cod} f= \mathrm{dom} g$, a composite arrow $g\circ f\colon \mathrm{dom} f \to \mathrm{cod} g$, satisfying the following associative law: for any $f\colon A\to B$, $g\colon B\to C$ and $h\colon C\to D$ , \[ h\circ (g \circ f) = (h\circ g) \circ f\: ;\]
    \item for each object $A$, an identity arrow $\id_A\colon A\to A$ satisfying the identity law: for any $ f \colon A\to B,\, g\colon{}B\to A $,
\[ \quad f \circ \id_A = f \text{ and } \id_A\circ g=g. \]
	\end{enumerate}
    \label{def:cat}
\end{defn}

\begin{defn}
    Given two categories $\mathcal C,\mathcal  D$ a \textbf{functor} $F\colon{}\mathcal C\to\mathcal  D$ is a map taking each $\mathcal C$-object $A$ to a $\mathcal D$-object $F(A)$ and each $\mathcal C$-arrow $f\colon A\to B$ to a $\mathcal D$-arrow $F(f) \colon F(A) \to F(B)$, such that
    \begin{enumerate}
        \item $F(g\circ f)=F(g)\circ F(f)$, whenever $g\circ f$ is well-defined;
        \item  $F(\id_A)=\id_{F(A)}$ for each object $A\in\mathcal  C$.
    \end{enumerate}
\end{defn}

The categories we are interested in possess an additional important structure: as well as sequential composition, as defined in Definition~\ref{def:cat}, they also possess a notion of \textit{parallel} composition. This is formally captured in the notion of a \textit{monoidal} category.

\begin{defn}
    A \textbf{monoidal category} is a category $\mathcal C$ equipped with the following:
    \begin{enumerate}
        \item $\mathcal C$ is equipped with a functor $\otimes\colon{} \mathcal C\times\mathcal  C \to \mathcal C$;
        \item $\mathcal C$ possesses a distinguished object $1\in \mathcal C$;
        \item there are isomorphisms $\alpha_{A,B,C}\colon{}(A\otimes B)\otimes C \to A\otimes (B\otimes C)$ for all objects $A,B,C$ which are natural in $A,B,C$;
        \item for all objects $A\in\mathcal  C$ there are isomorphisms, the \textit{left-} and \textit{right-unitors}, $\lambda_A\colon{} 1\otimes A \to A$ and $\rho_A\colon{} A\otimes 1 \to A$, both natural in $A$;
        \item these morphisms are subject to the \textit{pentagon equation}:
        \[
\begin{tikzcd}
  & ((A \otimes B) \otimes C) \otimes D \arrow[dl, "\alpha_{A,B,C} \otimes \text{id}_D"'] \arrow[dr, "\alpha_{A\otimes B,C,D}"] & \\
  (A \otimes (B \otimes C)) \otimes D \arrow[d, "\alpha_{A,B\otimes C,D}"'] & & (A \otimes B) \otimes (C \otimes D) \arrow[d, "\alpha_{A,B,C\otimes D}"] \\
  A \otimes ((B \otimes C) \otimes D) \arrow[rr, "\text{id}_A \otimes \alpha_{B,C,D}"'] & & A \otimes (B \otimes (C \otimes D))
\end{tikzcd}
\]

        \item and the \textit{triangle equation}:
        \[
        \begin{tikzcd}
             (A\otimes 1)\otimes B \arrow[rr, "\alpha_{A,1,B}"] \arrow[dr, "\rho_A \otimes 1_B"]& & A\otimes (1\otimes B) \arrow[dl,"1_A \otimes \lambda_B"]\\
             &  A\otimes B &
        \end{tikzcd};
        \]
        \item the category is furthermore a \textbf{symmetric} monoidal category if it is equipped with natural isomorphisms $s_{A,B}\colon{}A\otimes B \to B\otimes A$ such that 
        $$s_{B,A}\circ s_{A,B}=\id_{A\otimes B}$$
        and the hexagon identity is satisfied:
                \[
\begin{tikzcd}
  (A \otimes B) \otimes C \arrow[d, "s_{A,B}\otimes\id_C"'] \arrow[r, "\alpha_{A,B,C}"]& A\otimes (B\otimes C) \arrow[r,"s_{A,B\otimes C}"]& (B \otimes C) \otimes A  \arrow[d, "\alpha_{B,C,A}"] \\
  (B\otimes A)\otimes C \arrow[r, "\alpha_{B,A,C}"']  & B\otimes(A\otimes C) \arrow[r,"\id_B\otimes s_{A,C}"']& B\otimes(C\otimes A)
\end{tikzcd}.
\] 
    \end{enumerate}
\end{defn}

Our key contribution will be describing how the BFFG algorithm \textit{preserves both sequential and parallel composition}. This will be formalised by proving that the BFFG algorithm gives rise to a functor, which in some cases is lax monoidal.
\begin{defn}\label{defn:laxmf}
    Given two monoidal categories $(\mathcal C,\otimes_\mathcal C, 1_\mathcal C, \alpha^\mathcal C, \lambda^\mathcal C, \rho^\mathcal C)$ and $(\mathcal D,\otimes_\mathcal D, 1_\mathcal D, \alpha^\mathcal D, \lambda^\mathcal D, \rho^\mathcal D)$, a \textit{lax monoidal functor} is a functor $F\colon{}\mathcal C\to\mathcal  D$, together with
    \begin{enumerate}
        \item a morphism $\epsilon\colon{}1_\mathcal D \to F(1_\mathcal C)$;
        \item for all objects $A,B\in\mathcal C$, morphisms $\mu_{A,B}\colon{}F(A)\otimes_\mathcal D F(B) \to F(A\otimes_\mathcal C B)$, natural in $A,B$:

\[
\begin{tikzcd}
    F(A)\otimes_\mathcal D F(B) \arrow[r, "\mu_{A,B}"] \arrow[d, "F(f) \otimes_{\mathcal D} F(g)"'] 
    &F(A \otimes_\mathcal C B)  \arrow[d, "F(f \otimes_{\mathcal C}  g)"] \\
    F(C) \otimes_{\mathcal D} F(D) \arrow[r, "\mu_{C,D}"] 
    & F(C  \otimes_{\mathcal C} D) 
\end{tikzcd};
\]

        \color{black}
        \item these satisfy associativity: for all objects $X,Y,Z\in \mathcal C$, the following diagram commutes:
        \[
\begin{tikzcd}[column sep = 14ex]
  (F(A) \otimes_\mathcal D F(B)) \otimes_\mathcal D F(C) \arrow[r, "\alpha^\mathcal D_{F(A),F(B),F(C)}"] \arrow[d, "\mu_{A,B} \otimes \text{id}"'] 
  & F(A) \otimes_\mathcal D (F(B) \otimes_\mathcal D F(C)) \arrow[d, "\mathrm{id}\otimes \mu_{ B,C}"] \\
    F(A \otimes_\mathcal C B) \otimes_\mathcal D F(C)  \arrow[d, "\mu_{A\otimes_\mathcal C B,C}"'] 
  & F(A) \otimes_\mathcal D F(B \otimes_\mathcal C C) \arrow[d, "\mu_{A, B \otimes_\mathcal C C}"] \\
  F((A \otimes_\mathcal C B) \otimes_\mathcal C C) \arrow[r, "F(\alpha^\mathcal C_{A,B,C})"'] 
  & F(A\otimes_\mathcal C ( B \otimes_\mathcal C  C))
\end{tikzcd}
\]

        \item and satisfy unitality: for each object $A\in\mathcal  C$, the following diagrams commute:
        \[
\begin{tikzcd}
    1_\mathcal D\otimes_\mathcal D F(A) \arrow[r, "\epsilon\otimes \mathrm{id}"] \arrow[d, "\lambda^\mathcal D_{F(A)}"'] 
    & F(1_\mathcal C)\otimes_\mathcal D F(A) \arrow[d, "\mu_{1_\mathcal C,A}"] \\
    F(A) 
    & F(1_\mathcal C \otimes A) \arrow[l, "F(\lambda^\mathcal C_A)"'] 
\end{tikzcd}
\]
    and
        \[
\begin{tikzcd}
    F(A)\otimes_\mathcal D 1_\mathcal D \arrow[r, "\mathrm{id}\otimes \epsilon"] \arrow[d, "\rho^\mathcal D_{F(A)}"'] 
    & F(A)\otimes_\mathcal D F(1_\mathcal C) \arrow[d, "\mu_{A,1_\mathcal C}"] \\
    F(A) 
    & F(A \otimes 1_\mathcal C) \arrow[l, "F(\rho^\mathcal C_A)"'] 
\end{tikzcd}.
\]
    \end{enumerate}
\end{defn}


\bibliography{literature_fixed}
\bibliographystyle{plainnat}

\end{document}